\documentclass[a4paper]{article}

\usepackage{amsmath}
\usepackage{amssymb}
\usepackage{amsfonts}
\usepackage{enumerate}
\usepackage{vmargin}
\usepackage[arrow,matrix]{xy}

\newcommand{\tens}{\otimes}
\newcommand{\mat}{\mathcal}
\newcommand{\N}{\ensuremath{\mathbb{N}}}
\newcommand{\K}{\ensuremath{\mathbb{K}}}
\newcommand{\M}{\ensuremath{\mathbb{M}}}
\newcommand{\Z}{\ensuremath{\mathbb{Z}}}
\newcommand{\B}{\ensuremath{\mathbb{B}}}

\newcommand{\C}{\ensuremath{\mathbb{C}}}
\newcommand{\T}{\ensuremath{\mathbb{T}}}
\newcommand{\Id}{\ensuremath{\mathrm{Id}}}

\newcommand{\ran}{\ensuremath{\mathop{\rm Im\,}}}

\newcommand{\sk}{\smallskip}

\renewcommand{\le}{\ensuremath{\leqslant}}
\renewcommand{\ge}{\ensuremath{\geqslant}}
\renewcommand{\leq}{\ensuremath{\leqslant}}
\renewcommand{\geq}{\ensuremath{\geqslant}}
\newcommand{\n}{\noindent}
\newcommand{\la}{\langle}
\newcommand{\ra}{\rangle}
\newcommand{\qed}{\hfill \vrule height6pt  width6pt depth0pt}

\newcommand{\F}{\mathcal F}
\newcommand{\pl}[3]{\mathring{#1}_{#2_1}\otimes\cdots\otimes
\mathring{#1}_{#2_{#3}}}
\newcommand{\plB}[3]{\mathring{#1}_{#2_1}\otimes_B\cdots\otimes_B
\mathring{#1}_{#2_{#3}}}
\newcommand{\sd}[1]{\bigoplus\limits_{\substack{#1\ge1 \\
i_1\not=\,i_2\not=\cdots \not=\,i_{#1}}}}
\newcommand{\sdd}[1]{\bigoplus\limits_{i_1\not=\,i_2\not=\cdots \not=\,i_{#1}}}
\newcommand{\inn}[2]{#1_1\not=\,#1_2\not=\cdots \not=\,#1_{#2}}
\newcommand{\m}{\mathring}

\newtheorem{thm}{Theorem}[section]

\newtheorem{prop}[thm]{Proposition}
\newtheorem{cor}[thm]{Corollary}
\newtheorem{lemma}[thm]{Lemma}
\newtheorem{fact}[thm]{Fact}
\newenvironment{rk}{\refstepcounter{thm}\noindent%
{\bf Remark \arabic{section}.\arabic{thm}} \ }{
\smallskip

}

\newenvironment{pf}[1][]{\noindent {\it Proof #1} : }{\hbox{~}\qed
\smallskip
}

\title{Khintchine type inequalities for reduced free products\\
and Applications}
\date{}
\author{Eric Ricard and Quanhua Xu}

\begin{document}

\maketitle


\begin{abstract}
 We prove  Khintchine type inequalities for words of a fixed
length in a reduced free product of $C^*$-algebras (or von Neumann
algebras). These inequalities imply that the natural projection
from a reduced free product onto the subspace generated by the
words of a fixed length $d$ is completely bounded with norm
depending linearly on $d$. We then apply these results to various
approximation properties on reduced free products. As a first
application, we give a quick proof of Dykema's  theorem
on the stability of exactness under
the reduced free product for $C^*$-algebras. We next study the
stability of the completely contractive approximation property
(CCAP) under reduced free product. Our first result in this
direction is that a reduced free product of finite dimensional
$C^*$-algebras has the CCAP. The second one asserts that a von
Neumann reduced free product of injective von Neumann algebras has
the weak-$*$ CCAP.  In the case of group $C^*$-algebras, we show
that a free product of weakly amenable groups with constant 1 is
weakly amenable.
\end{abstract}

\makeatletter
\renewcommand{\@makefntext}[1]{#1}
\makeatother
\footnotetext{\noindent Laboratoire de Math\'{e}matiques,
 Universit\'{e} de Franche-Comt\'{e},
 25030 Besan\c con, cedex - France\\
 ricard@math.univ-fcomte.fr,\hskip 0.5cm qx@math.univ-fcomte.fr
 }

\begin{section}{Introduction and Background}
This paper deals with the reduced free product of $C^*$-algebras
(and of von Neumann algebras). The construction of reduced free
product was introduced independently by Avitzour \cite{Avit} and
Voiculescu \cite{Voi} (see also a previous work by Ching
\cite{ching}). Since then it has been considerably developed and
becomes today an independent direction of research, free
probability theory. This theory has many interactions with other
directions such that quantum probability, operator algebras and
operator spaces, and turns out to be an efficient tool notably for
the two last theories.

\medskip

Our starting point is Haagerup's inequality \cite{Haa}.
Let $\mathbb F_n$ be a free group on $n$ generators $g_1,...,g_n$
and $W_d$ be the subset of $\mathbb F_n$ of words of length $d$.
Then for any family $\{x_w\}$ of complex numbers
 $$\big( \sum_{w\in W_d} |x_w|^2\big)^{1/2}\le
 \big\|\sum_{w\in W_d } x_w \lambda(w)\big\|_{C_\lambda^*(\mathbb F_n)}
 \le (d+1)\big( \sum_{w\in W_d} |x_w|^2\big)^{1/2}\,,$$
where $C_\lambda^*(\mathbb F_n)$ is the reduced $C^*$-algebra of
$\mathbb F_n$ generated by the left regular representation of
$\mathbb F_n$ on $\ell_2(\mathbb F_n)$. The case of $d=1$ goes
back to Leinert \cite{leinert}. Regarding a family of free
generators, or more generally, $\{\lambda(w)\}_{w\in W_d}$, as a
lacunary set, the inequality above can be interpreted as a
Khintchine type inequality in $L_\infty$. Note that such a
phenomenon cannot occur in the commutative setting, namely, there
does not exist any infinite lacunary sequence in an abelian group
which generates in $L_\infty$ a subspace isomorphic to $\ell_2$.

Leinert's inequality was extended to the case of operator valued
coefficients by Haagerup and Pisier \cite{HP} :
 $$\big\|\sum_{i=1}^n \lambda(g_i) \tens x_i
 \big\|_{C_\lambda^*(\mathbb F_n)\tens \M_m}
 \le 2\max \Big\{\big\|\sum_{i=1}^n x_i^*x_i\big\|_{\M_m}^{1/2}
 \; ,\;\big\|\sum_{i=1}^n x_ix_i^*\big\|_{\M_m}^{1/2}\Big\},$$
where $x_i\in \M_m$ and $\M_m$ stands for the algebra of $m\times
m$ complex matrices. The converse inequality (with constant 1) is
easy. Later, Buchholz \cite{Buc} found a right formulation for the
operator valued version of Haagerup's inequality.

The reduced $C^*$-algebra of a free group is an important example
of reduced free product algebras. Thus it is natural to attempt to
transfer the previous inequalities to reduced free product. This
was done in the case of length 1 by Voiculescu \cite{V2} for the
scalar-valued case and by Junge \cite{ju} for the vector-valued
(or amalgamated) case. More precisely, we have the following free
product version of Haagerup-Pisier's inequality (without
amalgamation). Let $(A_i,\phi_i)_{i\in I}$ be a family of
$C^*$-algebras equipped with states $\phi_i$ whose GNS
constructions are faithful. Let $\mat A=*_{i\in I}(A_i,\phi_i)$ be
the associated reduced free product. Then for $a_i\in
(A_i,\phi_i)$ with $\phi_i(a_i)=0$ and $x_i\in\M_m$, we have
 \begin{eqnarray*}
 \big\|\sum_i a_i\tens x_i\big\|_{\mat A\tens \M_m}
 &\le& 3 \max\,\Big\{ \max_i\|a_i\|_{A_i}\,\|x_i\|_{\M_m}\;,\,\\
 & & \qquad \big\|\sum_i \phi_i(a_i^*a_i) x_i^*x_i
 \big\|_{\M_m}^{1/2}\;,\; \big\|\sum_i \phi_i(a_ia_i^*) x_ix_i^*
 \big\|_{\M_m}^{1/2}\,\Big\}.
 \end{eqnarray*}
Again the converse inequality is easy to be checked.

One of the main results of this paper is the extension of the
inequality above to an arbitrary fixed length $d$, i.e. the free
product version of Buchholz's inequality. The relevant constant is
then $2d+1$. This allows to show that the subspace generated by
the words of length $d$ in a reduced free product is complemented
with a constant depending linearly on $d$. We should emphasize
that this linear (or polynomial) dependence on $d$ is crucial for
applications. These results will be proved in the following two
sections. We will use the formalism of operator space theory,
notably, the part concerning row, column Hilbertian spaces and
Haagerup tensor product.

Section~4 presents some applications of the results just
mentioned. The common topic is approximation property in various
senses. For instance, the previous Khintchine type inequalities
can be used to provide a simple operator theoretic proof of the
stability of exactness under reduced amalgamated free product, a result due to
Dykema \cite{D} (another proof  was given Dykema-Shlyakhentko
\cite{DS} using Cuntz-Pimser algebra).

The main motivation of the paper comes, however,  from the open
problem whether the completely contractive approximation property
(CCAP) is preserved by reduced free product. Although we cannot
completely solve this problem,  we do provide positive solutions
in some particular cases.  Concerning the von Neumann algebra
reduced free product, we get a rather satisfactory solution: a von
Neumann reduced free product of injective von Neumann algebras
with respect to any normal states has the weak$^*$-CCAP.

On the other hand, in the case of group algebras, we completely
solve the problem above: a free product of weakly amenable
discrete groups with constant 1 is still weakly amenable with
constant 1. This last result is an improvement of a previous
theorem by Bo{\.z}ejko and Picardello \cite{BP}, which asserts
that a free product of amenable discrete groups is weakly
amenable.

 From the operator space point of view, the previous problem seems
quite natural. Around the same topic, a result of Sinclair and
Smith \cite{SiSm} states that the CBAP for $C^*$-algebras is stable
under crossed product with discrete amenable groups. More
recently, Dykema and Smith \cite{DSm} proved that Cuntz-Pimsner
algebras constructed over $C^*$-algebras with the CBAP also have
the CBAP.

\medskip

In the rest of this section we briefly recall the construction of
reduced free product. We will use standard notations and notions in
the theory of free products and operator space theory. Our
references are \cite{VDN, ER2, Pis}.

Throughout the paper, $(A_i,\phi_i)_{i\in I}$ will be a family of
unital $C^*$-algebras with distinguished states $\phi_i$ whose GNS
constructions $(\pi_i, H_i, \xi_i)$ are faithful. The cardinality
of $I$ will be often denoted by $N$ (it can be an infinite
cardinal number). Recall that $H_i$ is the Hilbert space
$L_2(A_i,\phi_i)$. Let $\m H_i=\xi_i^\bot$. There are natural maps
$\hat{} :A_i\to H_i$ such that $\phi_i(a^*b)=\la\hat a ,\; \hat
b\ra_{H_i}$.  $H_i^{op}$ is the Hilbert space obtained from $A_i$
with respect to the sesquilinear form $(a, b)=\phi_i(ab^*)$. With
this definition, we have $H_i^{op*}=H_i$ via the following duality
 $$\la\hat a,\;\hat b\ra_{(H_i^{op},\; H_i)}=\phi_i(ab).$$

Let $\F$ be the Hilbert Fock space associated to the free product:
 $$\F=\mathbb{C}\cdot\Omega ~~\bigoplus~\sd n \pl H i n.$$
 This Hilbert space has a natural gradation given by the direct sum.
It would be helpful to think elementary tensors of the form
$h_1\tens\cdots\tens h_n \in\pl H i n$ as words of length $n$ in
letters coming from $\m H_i$'s; being the empty word, $\Omega$ has
length 0. The set of all words is a linearly dense subset of $\F$.

Following Voiculescu \cite{VDN}, each algebra $A_i$ acts non
degenerately on $\F$ from the left. More precisely, for $a\in \m
A_i$ and $h_1\tens\cdots\tens h_n \in \pl H i n$ a word of length
$n\ge 0$
 $$ a \cdot ( h_1\tens\cdots \tens h_n)= \left\{\begin{array}{ll}
 \hat a \tens h_1\tens\cdots \tens h_n & \textrm{ if } i\neq i_1 \\
 \big[a\cdot h_1 - \la\xi_i,\, a\cdot h_1\ra\xi_i\big]
 \tens h_2\cdots \tens h_n   &\textrm{ if } i= i_1 \\
 \hskip 0.5cm+\la\xi_i,\, a\cdot h_1\ra h_2 \tens \dots  \tens h_n &
 \end{array}\right.$$
Thus, the action of $a$ on a word $w$ can be seen as divided into
three parts : $a$ adds a letter $\hat a$ at the left of $w$, we
will say that it is a creation, or $a$ acts on the first letter of
$w$, we will call it a diagonal action, or $a$ removes the first
letter of $w$, this is an annihilation. This terminology is
consistent with the gradation of $\F$;  somehow, $a$ is just a
tri-diagonal block operator.

The faithfulness of the GNS constructions of the $(A_i,\phi_i)$
ensures that the above representations  of $A_i$ on $\F$ are
faithful. Thus there is a copy of the algebraic free product
 $$ A=\C 1 \bigoplus \sd d \pl A i d$$
in $\B(\F)$, the algebra of all bounded operators on $\F$.  The
reduced free product of the family $(A_i, \phi_i)_{i\in I}$ is the
$C^*$-algebra generated by these actions. It is just the closure
of $A$ in $\B(\F)$. For convenience, it will be denoted by
 $$\mathcal A=*_{i\in I} (A_i,\phi_i).$$
The states $\phi_i$ determine a state $\phi$ on $\mathcal A$ given
by:
 $$\phi(1)=1 \qquad \textrm{and}\qquad \phi(a_1\tens\cdots\tens a_d)
 =0$$
for $d\geq 1$, and $a_1\tens\cdots\tens a_d\in\pl A i d$ with $\inn i d$. As
usual, each $A_i$ is naturally considered as a subalgebra of
$\mathcal A$. Then the restriction of $\phi$ to $A_i$ coincides
with $\phi_i$.  Recall that the family $\{A_i\}_{i\in I}$ is free
in $(\mat A, \phi)$.

 As $\F$, the algebraic free algebra $A$ is naturally graded. We will
denote its homogeneous part of degree $d$ by $\Sigma_d$:
 $$\Sigma_d =\bigoplus\limits_{\inn i d}\pl A i d.$$
The completion of $\Sigma_d$ in $\mathcal A$ will be called
$\mathcal A_d$:
 $$ \mathcal A_d= \overline{\Sigma_d}^{\mathcal A}.$$

As for elements in the free Fock space, we will often refer to
elements in $\m A_i$ as letters and elementary tensors in
$\Sigma_d$ as words. Thus viewed in $\mat A$, a word
$a_1\tens\cdots\tens a_d$ is also equal to the product $a_1\cdots
a_d$. On the other hand, $\mathcal A_d$ is the closed subspace of
$\mathcal A$ generated by the words of length $d$. We will call it
the homogeneous subspace of degree $d$.

\smallskip

The construction above can also be done in the category of von
Neumann algebras. Let $(M_i,\phi_i)$ be von Neumann algebras with
distinguished normal states whose GNS constructions are
faithful. Then the von Neumann algebra reduced free product of the
$(M_i,\phi_i)$ is the weak-* closure of $ *_i (M_i,\phi_i)$ in
$\B(\mathcal F)$, which will be denoted by $(\overline{\mat M},
\phi)=\overline *_i (M_i,\phi_i)$. Again, the $M_i$ are regarded
as von Neumann subalgebras of $\overline{\mat M}$, and then the restriction
of $\phi$ to $M_i$ is equal to $\phi_i$. The homogeneous subspace
of degree $d$ of $\overline{\mat M}$ is the weak*-closure of
$\Sigma_d$ above (with $A_i$ replaced by $M_i$). It will be
denoted by $\overline{\mat M}_d$.

\medskip

Concerning operator spaces, we will need only very few notions
beyond the basic definitions (completely bounded maps, minimal
tensor norm). If $H$ is a Hilbert space, we use the notation $H_C$
for the column operator space structure on $H$, that is obtained
by the obvious inclusion $H\subset \B(\C,H)$, the space of bounded
maps from $\C$ to $H$. Its row counterpart, $H_R$ comes from the
inclusion $H\subset \B(H^*,\C)$. For other unexplained definitions
(e.g. Haagerup tensor product), we refer to \cite{ER2, Pis}.
$\M_n$ stands for the full algebra of $n\times n$ complex
matrices and $\K$ for compact operators.

\medskip

In the remainder of the paper, all notations just introduced will
be kept to have the previous meanings, unless explicitly stated
otherwise.

\end{section}

\begin{section}{Khintchine Inequalities}

This section is devoted to the Khintchine type inequalities for
$\Sigma_d$. The following first lemmas are well known and will be
the basic building blocks. The underlying idea is quite simple, it
consists in decomposing the operators in $\Sigma_d$ in an
appropriate way with respect to the gradation of $\F$. As in
section~1, $(A_i, \phi_i)_{i\in I}$ denotes a family of
$C^*$-algebras with distinguished states $\phi_i$ whose GNS
constructions $(\pi_i, H_i, \xi_i)$ are faithful. $(\mathcal A,\,
\phi)$ is the associated reduced free product.

\smallskip

For each $k\in I$, let $P_k$ be the projection from $\F$ onto the
subspace
 $$\F_k=\bigoplus\limits_{\substack{n\ge 1 \\ k=\inn i n}} \pl H in$$
and $P^\bot_k$ its complement. Recall that $\Sigma_1$ is just the
direct sum $\oplus_i \m A_i$. Let $\sigma: \Sigma_1\to \B(\F)$ be
defined by $\sigma(a)=P_k a P_k$ if $a\in \m A_k$.
Actually, $\sigma$ can be defined from $\ell_\infty((A_i))$ to
$\B(\F)$ by the same formula.

\begin{lemma}\label{linfi}
The map $\sigma$ extends to a complete contraction from
$\ell_\infty( ( A_i))$ to $\B(\F)$.
\end{lemma}

\begin{pf}
It suffices to prove that for $a_{k,i}\in A_k$ and $m_{k,i}\in \M_d$
($d\in\N$)
 $$\big\|\sum_{k,i} P_ka_{k,i}P_k\tens m_{k,i}\big\|_{\B(\F)\tens \M_d}
 \le\sup_{k} \big\|\sum_i a_{k,i}\tens m_{k,i}\big\|_{\B(H_k)\tens \M_d}.$$
Since the $P_k$ are mutually orthogonal, we have
 \begin{eqnarray*}
 \big\|\sum_{k,i} P_ka_{k,i}P_k\tens m_{k,i}\big\|_{\B(\F)\tens \M_d}
 &=&\sup_k \big\| \sum_i P_ka_{k,i}P_k\tens m_{k,i}\big\|_{\B(\F)\tens \M_d}
 \\
 &\leq & \sup_k \big\|\sum_i a_{k,i}\tens m_{k,i}\big\|_{\B(\F)\tens
 \M_d}\\&=& \sup_{k} \big\|\sum_i a_{k,i}\tens m_{k,i}\big\|_{A_k\tens \M_d}.
 \end{eqnarray*}
The last equality occurs as the embeding $A_k\subset B(\F)$ is a complete
isometry.
\end{pf}

\begin{lemma}\label{annul}
 Let $a\in A_k$. Then $P^\bot_kaP^\bot_k=\phi_k(a) P^\bot_k$.
\end{lemma}

\begin{pf}
We can assume $\phi(a)=0$ (and so $a\in\m A_k$). The range of
$P^\bot_k$ is the span of elementary tensors
$h=h_{1}\tens\cdots\tens h_{n}$, where $h_1$ does not belong to
$\m H_k$. However, for such tensors, $a\cdot h\in \mathcal F_k$,
and so $P^\bot_k(a\cdot h)=0$.
\end{pf}

Let $L_1$ be the operator space in $\B(\F)$ spanned by $(P_k \m A_k
P^\bot_k)_{k\in I}$.

\begin{lemma}\label{Colonne}
 We have a complete isometry
$$L_1\approx \big(\bigoplus_{k\in I} \m H_k\big)_C.$$
 More precisely, for $a_{k,i}\in \m A_k$ and $m_{k,i} \in \M_d$:
 $$\big\| \sum_{k,i} P_k a_{k,i} P^\bot_k \tens m_{k,i}\big\| =
 \big\| \sum_{k} \sum_{i,j}\phi_k(a_{k,i}^* a_{k,j})m_{k,i}^* m_{k,j}
\big\|^{1/2}.$$
Moreover, the natural map $\theta_1:\mathcal A_1\to L_1$ defined
by $\theta_1(a)=P_kaP^\bot_k$ if $a\in \m A_k$ is a complete
contraction.

\end{lemma}
\begin{pf}
Let $a_{k,i}\in\m A_k$. Since the $P_k$ are mutually orthogonal
projections, thanks to the previous lemma, we have:
\begin{eqnarray*}
&&\Big(\sum_{k,i}  P_k a_{k,i} P^\bot_k \tens m_{k,i} \Big)^*
\Big(\sum_{k,i} P_k a_{k,i} P^\bot_k \tens m_{k,i}\Big)
\\&=& \sum_{k} \Big( \sum_i P^\bot_ka_{k,i}^*\tens m_{k,i}^*\Big)
(P_k \tens \Id)
\Big( \sum_i a_{k,i} P^\bot_k \tens m_{k,i}\Big)\\
&\leq& \sum_{k} \Big( \sum_i P^\bot_ka_{k,i}^*\tens m_{k,i}^*\Big)
\Big( \sum_i a_{k,i} P^\bot_k \tens m_{k,i}\Big)
\\&=& \sum_{k} P^\bot_k\tens
\sum_{i,j} \phi_k(a_{k,i}^*a_{k,j}) m_{k,i}^*m_{k,j}\\&
\leq&\sum_{k}\sum_{i,j} \phi_k(a_{k,i}^*a_{k,j})\Id \tens
m_{k,i}^*m_{k,j}.
\end{eqnarray*}
The last inequality comes from the fact that $\sum_{i,j}
\phi_k(a_{k,i}^*a_{k,j}) \tens m_{k,i}^*m_{k,j}$ is a positive operator
for any $k$.

 This gives the majoration. For the minoration, one only needs to
determine the action of the first term above  on $\Omega$.

 The second assertion follows from the first one for the natural map
$\mathcal A_1 \to \big(\bigoplus_{k\in I} \m H_k\big)_C$ is a
complete contraction.
\end{pf}

Passing to adjoints, and letting $K_1=L_1^*$, we obtain similarly the

\begin{cor}\label{Ligne}
We have completely isometrically
 $$K_1\approx \big(\bigoplus_{k=1}^N \m H^{op}_k\big)_R.$$
With, for $a_{k,i}\in \m A_k$ and $m_{k,i} \in \M_d$ :
 $$\big\| \sum_{k,i} P^\bot_k a_{k,i} P_k \tens m_{k,i}\big\| =
 \big\| \sum_{k}\sum_{i,j} \phi_k(a_{k,i} a_{k,j}^*)m_{k,i} m_{k,j}^*
 \big\|^{1/2}.$$
Moreover, the natural map $\rho_1:\mathcal A_1\to K_1$ defined
by $\rho_1(a)=P^\bot_kaP_k$ if $a\in \m A_k$ is a complete contraction.
\end{cor}

Algebraically, we can identify $\Sigma_d$ with a subspace of
$ \Sigma_1^{\tens^d}$.
We introduce an operator space structure on $ \Sigma_1^{\tens^d}$ via the
following inclusion
$$\iota : \left\{\begin{array}{ccccc}
 \Sigma_1^{\tens^d} & \to & \bigoplus_{k=0}^d L_1^k\tens_h K_1^{d-k} &
\bigoplus_\infty &
\bigoplus_{k=0}^{d-1}L_1^k\tens_h\ell_\infty((A_i)) \tens_h K_1^{d-k-1}\\
a& \mapsto& \Big((\theta_1^k\tens\rho_1^{d-k}(a))_{k=0}^d&, &
(\theta_1^k\tens \Id\tens \rho_1^{d-k-1}(a))_{k=0}^{d-1}\Big)
\end{array} \right.,$$
where $L_1^k=L_1^{\tens_h^ k}$ and $K_1^k=K_1^{\tens_h^ k}$. Here
the direct sum is in the $\ell_\infty$-sense.\\
 The big sum appearing on the right  will be denoted by $X_d$ in the
sequel. The induced operator space structure obtained on
$\Sigma_d$ after completion is denoted by $E_d$ (the fact that it
is indeed a norm will follow from the following theorem). We
denote by $\kappa$ the inclusion from $\Sigma_d$ to $\mathcal
A_d$.

 We also need to introduce $\overline X_d$ as
 $$ \bigoplus_{k=0}^d \B(K_1^{*(d-k)},\; L_1^k) \bigoplus_\infty\,
 \bigoplus_{k=0}^{d-1}\B(K_1^{*(d-k-1)},\;L_1^k)\,\overline \tens\,
 \ell_\infty((A_i'')).$$
By virtue of elementary properties of  Haagerup tensor product,
there is a natural embedding of $X_d$ in $\overline X_d$.

The following is the Khintchine inequality for $\Sigma_d$. Note
that \cite{jpx} contains a variant of this inequality as well as
its generalization to noncommutative $L_p$-spaces.

\begin{thm}\label{khintchine}
We have a complete isomorphism $E_d\approx \mathcal A_d$. More
precisely, for any $n\ge 1$ and $x\in \M_n(\Sigma_d)$, we have :
 $$ \|\iota(x)\|_{\M_n(E_d)}\le \|\kappa(x)\|_{\M_n(\mathcal A_d)}
 \le(2d+1) \|\iota(x) \|_{\M_n(E_d)}.\qquad (K_d)$$
\end{thm}
\begin{pf}
 The proof will consist in
 constructing two maps $\Pi_d$ and $\Theta_d$
with $\|\Pi_d\|_{cb}\le 2d+1$, $\|\Theta_d\|_{cb}\le 1$ such that
 the following diagram commutes

$$\xymatrix{
   &   \mathcal A_d \ar@{^{(}->}[r] & \B(\mathcal F)\ar[rdd]^{\Theta_d}& \\
\Sigma_d \ar[ur]^\kappa \ar[dr]_\iota & & &\\
 & E_d  \ar@{^{(}->}[r] & X_d \ar@{^{(}->}[r]
\ar[uu]^{\Pi_d\;\;} & \overline X_d
}$$

With these identifications $\Sigma_d$ is a dense subspace of both
$\mathcal A_d$ and $E_d$, thus to get an isomorphism in the
theorem, it suffices to prove that the norms induced on $\Sigma_d$
are equivalent. This boils down to the norm estimates in $(K_d)$.

\medskip\n{\em Majoration}~:~
We start with the  upper estimate  and the definition of $\Pi_d$.
For any $0\le k\le d$, the product map
 $$L_1^k\tens_h  K_1^{d-k} \to \B(\F)$$
is completely contractive by the very definition of the Haagerup
tensor product. In the same way, the map
 $$\begin{array}{ccc}
 L_1^k\tens_h\ell_\infty(( A_i)) \tens_h K_1^{d-k-1}& \to &\B(\F) \\
 x_1\tens\cdots \tens x_d &\mapsto &
 x_1\cdots x_k\sigma(x_{k+1})x_{k+2}\cdots x_d
 \end{array}$$
is a complete contraction since $\sigma$ is. Hence, we can define
a map $\Pi_d : X_d \to \B(\F)$ as the formal sum of the previous
product maps. It is completely bounded with norm less than $2d+1$.
If we take $a_1\tens\cdots \tens a_d \in \pl A i d$, with $\inn i
d$, we have,
 \begin{eqnarray*}
 \Pi_d(\iota ( a_1\tens\cdots \tens a_d))
 &=&\sum_{k=0}^d P_{i_1} a_1P^\bot_{i_1}\,\cdots\,
 P_{i_k} a_kP^\bot_{i_k}\,P^\bot_{i_{k+1}} a_{k+1}P_{i_{k+1}}\,
 \cdots\, P^\bot_{i_{d}} a_{d}P_{i_{d}}\\
 &+&\sum_{k=0}^{d-1} P_{i_1} a_1P^\bot_{i_1}\,\cdots\,
 P_{i_{k}} a_kP^\bot_{i_k}\,P_{i_{k+1}} a_{k+1}P_{i_{k+1}}\,
 P^\bot_{i_{k+2}} a_{k+2}P_{i_{k+2}}\,\cdots P^\bot_{i_{d}} a_{d}P_{i_{d}}.
 \end{eqnarray*}

 To get the upper estimate, it suffices to prove that the above
expression is exactly $ a_1\tens\cdots \tens a_d$ viewed in the
free product.
\begin{fact}\label{factdec}
In $\B(\F)$, we have the identities
 \begin{eqnarray*}
  \kappa(a_1\tens\cdots \tens a_d)&=& a_1 \,a_2\;...\;a_d\\
  &=& \prod_{k=1}^d (P_{i_k}+P^\bot_{i_k}) a_k
 (P_{i_k}+P^\bot_{i_k}) \\
 &=& \sum_{k=0}^d P_{i_1} a_1P^\bot_{i_1}\,\cdots \,
 P_{i_k} a_kP^\bot_{i_k}\, P^\bot_{i_{k+1}} a_{k+1}P_{i_{k+1}}\,
 \cdots \,P^\bot_{i_{d}} a_{d}P_{i_{d}}\\
 &+&
 \sum_{k=0}^{d-1} P_{i_1} a_1P^\bot_{i_1}\,\cdots\,
 P_{i_{k}} a_kP^\bot_{i_k}\,P_{i_{k+1}} a_{k+1}P_{i_{k+1}}\,
 P^\bot_{i_{k+2}} a_{k+2}P_{i_{k+2}}\,\cdots\, P^\bot_{i_{d}} a_{d}P_{i_{d}}.
 \end{eqnarray*}
 \end{fact}

To prove this fact we need to show  that when we expand the
product $a_1 \,a_2\;...\;a_d$, all terms vanish but those
appearing in the definition of $\Pi_d$.

There are $4^d$ terms in the development. Each term (called also a
word below) is a product of $d$ factors of the form
$Q_{i_j}a_jQ'_{i_j}$ with $Q_{i_j},
Q'_{i_j}\in\{P_{i_j},\;P^\bot_{i_j}\}$.  First, due to Lemma
\ref{annul}, all terms containing a factor of the form
$P^\bot_{i_j}a_jP^\bot_{i_j}$ are $0$. We analyze the terms in the
development that can contribute to the sum, i.e. those that do not
contain a factor of the form $P^\bot_.a_.P^\bot_.$ and show that
they are exactly the $2d+1$ terms from $\Pi_d$.

We proceed case by case, corresponding to the position of a letter
$P^\bot_.$.

 \begin{enumerate}[1)]
\item
 Assume that there is a $P^\bot_.$ immediately at the left of an $a_.$
in a given term. Then there must be a $P_.$ at the right of this
$a_.$.
 Consequently, the next factor on the right
must be of the form $P^\bot_.a_.P_.$ as the $P'$s are mutually
orthogonal. So the word is of the form $\cdots
P^\bot_.a_.P_.P^\bot_.a_.P_. \cdots P^\bot_.a_.P_.$.

Reading from left to right, consider now the first $P^\bot_.$
immediately at the left of an $a_.$, say, it is at position $k$.
Since it is the first with that property, for the factor at
position $(k-1)$, we have only two possibilities:
\begin{itemize}
 \item it is $P_{i_{k-1}} a_{k-1} P^\bot_{i_{k-1}}$, then
the whole word must be
 $$P_{i_1} a_1P^\bot_{i_1}\,\cdots\, P_{i_{k-1}}a_{k-1}P^\bot_{i_{k-1}}\,
 P^\bot_{i_{k}} a_{k}P_{i_{k}}\,\cdots\,
 P^\bot_{i_{d}}a_{d}P_{i_{d}}$$
according to the previous observations.
 \item it is $P_{i_{k-1}} a_{k-1} P_{i_{k-1}}$, then the letter at
the right of $a_{k-2}$ must be a $P^\bot$, so the whole word is
 $$P_{i_1} a_1P^\bot_{i_1}\,\cdots\,
 P_{i_{k-2}}a_{k-2}P^\bot_{i_{k-2}}\,
 P_{i_{k-1}} a_{k-1}P_{i_{k-1}}\,P^\bot_{i_{k}} a_{k}P_{i_{k}}\,
\cdots\, P^\bot_{i_{d}} a_{d}P_{i_{d}}.$$
\end{itemize}

\item Assume that we are not in the previous situation but
there is a $P^\bot_.$ immediately at the right
 of an $a_.$. Then there must be a $P_.$ at the left of this $a_.$.
So at the right of the previous $a_.$ there must be a $P^\bot_.$.
Thus the word is of the form $P_.a_.P^\bot_.\,P_.a_.P^\bot_.\,
\cdots\, P_.a_.P^\bot_.\, \cdots $.

We can consider the last $P^\bot_.$ at the right of an $a_.$.
Arguing as above leads to only two possibilities. Compared to the
development of $\Pi_d(\iota(a_1\tens \cdots a_d))$, they
correspond to the terms $k=d$ in the first sum and $k=d-1$ in the
second.

\item If there is no $P^\bot_.$ at all, then since the $P_.$ are
mutually orthogonal, we must
have $d=1$. For this length, the result is obvious.
\end{enumerate}

Therefore, we have proved the fact and thus the majoration in
$(K_d)$.

\medskip\n{\em Minoration}~:~  Now we turn to the lower estimate in $(K_d)$.
First, we fix $0\le k\le d$ and
prove that for $a\in \M_n(\Sigma_d)$, we have
$$\| a\|_{\M_n(\mathcal A_d)}\ge \| \theta_1^k\tens\rho_1^{d-k}(a)\|
_{\M_n(L_1^k\tens_hK_1^{{d-k}})}.$$
Since $L_1$ is a column operator space and $K_1$ a row operator space, we have
as an operator space
 $$L_1^k\tens_hK_1^{d-k} \subset \B({K_1^{d-k}}^*,\; L_1^k)\quad
 (\mbox{completely isometrically}).$$
By the obvious identification ${K_1^{d-k}}^*=L_1^{d-k}$,
 $L_1^k\tens_hK_1^{d-k}$ is so identified with a subspace of
 $\B((\bigoplus\limits_{i\in I} \m H_i)^{\tens^{d-k}},
 \;(\bigoplus\limits_{i\in I} \m H_i)^{\tens^k})$.

 For elementary tensors
$a=a_1\tens\cdots \tens a_d\in \pl A i d \in \Sigma_d$ with $\inn
i d$ and $\hat b_{1}\tens\cdots \tens \hat b_{d-k}\in \pl H j
{d-k}$ with $\inn j {{d-k}}$, we have the formula
 $$\theta_1^k\tens\rho_1^{d-k}(a)(\hat b_{1}\tens\cdots \tens\hat b_{d-k})
 =\phi(a_db_{1})\cdots \phi(a_{k+1}b_{d-k})\,
 \hat a_1\tens\cdots  \tens \hat a_k. \qquad  (*)$$

The operators we are interested in are those from
$\theta_1^k\tens\rho_1^{d-k}(\Sigma_d)$.  The conditions on the
indices in the sum  in $\Sigma_d$ imply that they do not act (i.e.
they vanish) on $(\sdd {d-k} \pl H i {d-k})^\bot$,  the complement
being taken in $(\bigoplus\limits_{i\in I} \m H_i)^{\tens^{d-k}}$.
Moreover, their ranges are contained in $\sdd {k} \pl H i {k}$.
Thus as an operator space, $\theta_1^k\tens\rho_1^{d-k}(\Sigma_d)$
is completely isometrically included in
 $$\B\big(\sdd {d-k} \pl H i {d-k},\;\sdd {k} \pl H i {k}\big).$$
Let $\mat P_n$ be the projection from $\F$ onto the subspace
generated by the words of length $n$. Then we claim that for any
$a\in\Sigma_d$
 $$ \mat P_{k} a_{|\ran(\mat P_{d-k})}
 = \theta_1^k\tens\rho_1^{d-k}(a).$$
This is easy by a length argument. By linearity, it suffices to
consider the case where $a=a_1\tens\cdots \tens a_d$ is an
elementary tensor. According to the decomposition in Fact
\ref{factdec}, $a$ acts on $\F$ as follows: either $a$ first
annihilates $q$ times and then creates $(d-q)$ times, or $a$ first
annihilates $q$ times, then acts once diagonally and finally
creates $(d-q-1)$ times.  It is clear that to pass from a word in
$\F$ of length $d-k$ to a word of length $k$, the latter case
cannot occur, and that in the former, $q$ must be equal to $d-k$,
i.e. $a$ must first annihilate $d-k$ times and then create $k$
times. This is exactly the formula $(*)$.

\sk

For the second kind of terms in the minoration, we apply the same
identification procedure. First, by standard results on Haagerup
tensor product, the operator space $L_1^k\tens_h\ell_\infty((A_i))
\tens_h K_1^{d-k-1}$ is naturally embedded in
$\ell_\infty((\B((\bigoplus\limits_{i\in I} \m
H_i)^{\tens^{d-k-1}},\; (\bigoplus\limits_{i\in I} \m
H_i)^{\tens^k})\tens_{\min} A_j)_j)$. So its norm is the supremum
of $N$ norms ($N$ being the cardinal of $I$). Let us fix $j$ and
concentrate on the norm for this $j$. Consider an element $\alpha
\in \M_n(\Sigma_d)$, and denote by $C_j(\alpha)$ its part whose
$(k+1)^{\rm th}$ letters are in $\m A_j$. Then the norm
corresponding to $j$ is
$$\| C_j(\alpha)\|_{\B(\sdd {d-k-1} \pl H i {d-k-1},\;\sdd {k} \pl H i {k}
)\tens_{\min} A_j \tens \M_n}.$$
Since in every word of
$C_j(\alpha)$, the two letters immediately before and after the
$(k+1)^{\rm th}$ belong  respectively to $\m A_i$ and $\m A_{i'}$
with $i\neq j$, $i'\neq j$, the norm of $C_j(\alpha)$ is the same
as the norm of a matrix of operators
 $$\bigoplus\limits_{i_1\not=\,i_2\not=\cdots \not=\,i_{d-k-1}\neq j}
 \pl H i {d-k-1} \tens  H_j\ \to\bigoplus
 \limits_{ i_1\not=\,i_2\not=\cdots \not=\,i_{k}\neq j}
 \pl H i {k} \tens  H_j.$$

After these preliminary observations, we now restrict our
attention to a fixed elementary tensor $\alpha=a=a_1\tens\cdots
\tens a_d\in \pl A l d \in \Sigma_d$ with $\inn l d$. Let
$p=l_{k+1}$. For $h=\hat b_{1}\tens\cdots \tens \hat
b_{d-k-1}\tens w \in \pl H j {d-k-1} \tens H_{j}$ with $\inn j{d-k-1}$,
 we have
 $$C_j((\theta_1^k\tens \Id\tens \rho_1^{d-k-1}(a)))h
 =\delta_{p,j}\,\phi(a_db_{1})\cdots \phi(a_{k+2}b_{d-k-1})\,
 \hat a_1 \tens\cdots \tens \hat a_k \tens
 \widehat{ a_{k+1}w}. \quad (**)$$

As previously, the operators of this type can be recovered
directly from $a$ using restrictions and compressions. To that
end, consider the two subspaces of $\F$ defined by
 $$S= \bigoplus\limits_{i_1\not=\,i_2\not=\cdots \not=\,i_{d-k-1}\neq j}
 \pl H i {d-k-1}\,\bigoplus\,  \bigoplus
 \limits_{i_1\not=\,i_2\not=\cdots \not=\,i_{d-k-1}\neq j}
 \pl H i {d-k-1}\tens \m H_j$$
 $$T=\bigoplus\limits_{i_1\not=\,i_2\not=\cdots \not=\,i_{k}\neq j}
 \pl H i {k} \tens \m H_j\,\bigoplus\, \bigoplus
 \limits_{i_1\not=\,i_2\not=\cdots \not=\,i_{k}\neq j} \pl H i {k}.$$
We have the following obvious identifications:
 $$S\approx \bigoplus
 \limits_{i_1\not=\,i_2\not=\cdots \not=\,i_{d-k-1}\neq j}
 \pl H i {d-k-1} \tens  H_j$$
 $$ T\approx  \bigoplus
 \limits_{i_1\not=\,i_2\not=\cdots \not=\,i_{k}\neq j}
 \pl H i {k} \tens  H_j.$$
To conclude, we just need to check that
 $$C_j((\theta_1^k\tens \Id\tens \rho_1^{d-k-1}(a))) = U a_{|S}\,,$$
where  $U$ is the projection from $\mathcal F$ to $T$.
 Let $h=\hat b_{1}\tens\cdots \hat b_{d-k-1}\tens w \in \pl H j {d-k-1}
\tens H_{j}$ as above, and let us determine the actions of both
operators on it.

The first remark is that both $C_j((\theta_1^k\tens \Id\tens
\rho_1^{d-k-1}(a)))h$ and $Ua\cdot h$ can be non zero only if
$b_1\in \m H_{l_d}, \cdots , b_{d-k-1}\in \m H_{l_{k+2}}$. For
$Ua\cdot h$, this follows by a length argument: otherwise, $a\cdot
h$ is a sum of words of length at least $k+2$. For the second
operator $C_j((\theta_1^k\tens \Id\tens \rho_1^{d-k-1}(a)))$, this
is clear by $(**)$. Now, we distinguish two cases with respect to
the value of $w$ :

\begin{itemize}
\item Assume $w=\xi_j$. Then viewed in $S$, $h=\hat b_{1}\tens\cdots \tens
 \hat b_{d-k-1}$. It then follows that
 $$a\cdot h=\phi(a_d b_{1})\cdots\phi(a_{k+2}b_{d-k-1})\,
 \hat a_1 \tens\cdots\tens \hat a_k\tens \hat a_p.$$
 Thus, this yields
 $$Ua\cdot h= \delta_{p,j}\,\phi(a_db_{1})\cdots\phi(a_{k+2}b_{d-k-1})\,
 \hat a_1 \tens\cdots \tens \hat a_k\tens\hat a_p.$$
By $(**)$, this is exactly
 $C_j((\theta_1^k\tens \Id\tens \rho_1^{d-k-1}(a)))h$.
\item Assume $w\in \m H_j$.
If $p\neq l_{k+1}$, then $a\cdot h$ is a finite sum of words of
length greater than $k+1$, so to prove the announced equality we
can assume that $p=l_{k+1}$. Then as above we easily recover
$(**)$ using the identification of $T$ with $\bigoplus\limits_{
i_1\not=\,i_2\not=\cdots \not=\,i_{k}\neq j} \pl H i {k} \tens
H_j$.
\end{itemize}
 This concludes the proof for the lower bound. If we sum up all
compressions (injections) and restrictions, we get the desired map
$\Theta_d:\B(\mathcal F)\to X_d$.
\end{pf}

\medskip

\begin{rk}\label{normalmaps}
We have chosen a fast way to construct the map $\Pi_d$ using the
Haagerup tensor product. Actually, it is possible to define an
extension $\Xi_d$ of $\Pi_d$ defined on $\overline X_d$ which is
$*$-weakly continuous. It suffices to notice that $\kappa(a)$ can
be recovered from $\iota(a)$ just using sums of
ampliations/restrictions/projections. We give a brief sketch
keeping the same notation as before.  For instance, let $t\in
\B(L_1^{d-k-1},\; L_1^k)\overline \tens A_j'' \subset
\B(L_1^{d-k-1}\tens H_j,\; L_1^k\tens H_j)$. Using some
compression, we can define an operator $m_{j,k}(t)\in \B(S,T)$ as
at the end of the proof above. Then tensorizing $m_{j,k}(t)$ with
the identity of $\mathcal F$ yields an operator $n_{j,k}(t)\in
\B(S\tens \mathcal F, T\tens \mathcal F)$. Let $u:\mathcal F \to
S\tens \mathcal F$ and $v: \mathcal F \to T\tens \mathcal F$ be
the natural partial isometries obtained from associativity of
tensor product. Put $l_{k,j}(t)=v^*n_{j,k}(t)u$. Summing over $j$
gives a completely contractive map from $\B(L_1^{d-k-1},\;
L_1^k)\overline \tens \ell_\infty(A_j'') \to \B(\mathcal F)$. We
can also perform the same kind of operations for elements in
$\B(L_1^{d-k},\; L_1^k)$. The addition of all those maps is the
normal map $\Xi_d$.
\end{rk}

\medskip

\begin{rk}
It is possible to write down a more concrete formula for the norm
in $E_d$. We do this only for $d=2$. A  typical element of
$M_n(\Sigma_2)$ is $x=\sum_{i_1\neq i_2} \sum_k m_{i_1,i_2}^k\tens
a_{i_1,k}^1\tens a_{i_2,k}^2$ (a finite sum) with
$m_{i_1,i_2}^k\in \M_n$ and $a_{i_j,k}^ {.}\in \m A_{i_j}$. The
norm of $a$ is then equivalent to the maximum of five terms
$$\begin{array}{c}
 \Big\| \sum_{i_1\neq i_2,k} \phi_{i_1}(a_{i_1,k}^{1*}a_{i_1,k}^1)
\phi_{i_2}(a_{i_2,k}^{2*}a_{i_2,k}^2) m_{i_1,i_2}^{k*}
m_{i_1,i_2}^k\Big\|_{\M_n}
 \\
 \Big\|\sum_k m_{i_1,i_2}^k\tens  a_{i_1,k}^1\tens
a_{i_2,k}^2\Big\|_{\M_n(L_1\tens_h K_2)}\\
 \Big\| \sum_{i_1\neq i_2,k} \phi_{i_1}(a_{i_1,k}^1a_{i_1,k}^{1*})
\phi_{i_2}(a_{i_2,k}^2a_{i_2,k}^{2*}) m_{i_1,i_2}^km_{i_1,i_2}^{k*}
\Big\|_{\M_n}
 \\
\sup_{i_1} \Big\|\sum_k\phi_{i_2}(a_{i_2,k}^2a_{i_2,k}^{2*}) m_{i_1,i_2}^k
m_{i_1,i_2}^{k*}
\tens a_{i_1,k}^1a_{i_1,k}^{1*}\Big\|_
{\M_n(A_{i_1})}\\
\sup_{i_2} \Big\|\sum_k\phi_{i_1}(a_{i_1,k}^{1*}a_{i_1,k}^1) \textmd{}
m_{i_1,i_2}^{k*}m_{i_1,i_2}^k\tens a_{i_2,k}^{1*}a_{i_2,k}^2\Big\|_
{\M_n(A_{i_2})}.
\end{array}$$

\end{rk}

\begin{rk}
When the  set $I$ is finite (with cardinal $N$), we can forget the
terms coming from $\oplus_{k=0}^d L_1^k\tens_h K_1^{d-k}$ in
$E_d$, but then we have to replace the constant $2d+1$ in $(K_d)$
by $(\sqrt N +1)d+1$. Indeed they are dominated by the other $d+1$
terms but we have to pay for the norm of the identity maps from
$\ell_\infty((A_i))$ to $L_1$ and to $K_1$.
\end{rk}

\medskip

\begin{rk}
 In the result of Buchholz \cite{Buc}, there are only $d+1$ terms in
his Khintchine inequality. This is not surprising, if one notices
that $\Sigma_d$ in the free product $\mathop *\limits_{1\le i\le
n}C(\mathbb T)$ and the span of $\lambda(w)$ for words of length
$d$ in $C^*_\lambda({\mathbb F}_n)$ are actually different. For
instance $\lambda(g_1^2g_2)$ has length 3 as a free word in
$C^*_\lambda(\mathbb F_n)$ but is in $\Sigma_2$ and hence has
length 2 according to our terminology. The diagonal actions mainly
explain that difference.
\end{rk}

\medskip
  We point out that the previous results
extend almost verbatim to  von Neumann algebras free product and
amalgamated free product. We conclude this section by a very brief
discussion on the former, and postpone the latter to the last
section.

\medskip

Let $(M_i,\phi_i)$ be von Neumann algebras with distinguished
normal states (with faithful GNS constructions). Then the von
Neumann reduced free product $\overline *_i (M_i,\phi_i)$ is the
weak-* closure of $ *_i (M_i,\phi_i)$ in $\B(\mathcal F)$. We
still keep the same notations as before, but we denote by
$\overline {\mathcal M}_d$ the weak-$*$ closure of $\Sigma_d$ in
$\B(\F)$.
 As $M_i''=M_i$,
$\overline X_d$, defined previously, is the weak-$*$ version of
$X_d$, that is
 $$ \bigoplus_{k=0}^d \B(K_1^{*(d-k)},\; L_1^k)\;\bigoplus_\infty\;
 \bigoplus_{k=0}^{d-1}\B(K_1^{*(d-k-1)},\; L_1^k)\,\overline \tens\,
 \ell_\infty((M_i)).$$
To be consistent with our previous approach, we can write it as
 $$\bigoplus_{k=0}^d L_1^k\tens_{eh} K_1^{d-k}
 \;\bigoplus_\infty\;
 \bigoplus_{k=0}^{d-1}L_1^k\tens_{eh}\ell_\infty((M_i))
 \tens_{eh}K_1^{d-k-1},$$
where $\tens_{eh}$ stands for the extended Haagerup tensor product
(see \cite{ER1,ER2}). Let $\overline E_d$ be the weak-$*$ closure
of $\Sigma_d$ in $\overline X_d$.

 It is straightforward to check that all maps ($\theta_1$,
$\rho_1$) considered earlier for $C^*$-algebras are weak-$*$
continuous in the von Neumann algebra setting.

Now the situation is very similar to that before and is summed up
by the diagram

 $$\xymatrix{
   &   \overline {\mathcal M}_d \ar@{^{(}->}[r] & \B(\mathcal F) \\
\Sigma_d \ar[ur]^\kappa \ar[dr]_\iota & & \\
 & \overline
E_d  \ar@{^{(}->}[r] & \overline X_d \ar[uu];[]_{\Xi_d\;\;}
\ar@<1ex>[uu]_{\;\; \Theta_d}}$$

The main difference with the $C^*$-algebra case comes from the
fact that the images of $\Sigma_d$ are only $*$-weakly dense. So
to conclude to any kind of isomorphism theorem, we should ensure
that the maps $\Theta_d$ and $\Xi_d$ are weak-$*$ continuous.  The
map $\Theta_d$ was obtained only from compressions/injections and
restrictions, so it is a normal map. By Remark \ref{normalmaps},
the map $\Xi_d$ defined earlier was also normal. Thus we have
obtained the

\begin{thm}\label{vnversion} The map
$j=\iota\kappa^{-1}$ defined on $\kappa(\Sigma_d)$ extends to  a
complete weak-$*$ continuous isomorphism between
 $\overline {\mat M}_d$ and $\overline E_d$.
More precisely, $\|j\|_{cb}\le 1$ and $\|j^{-1}\|_{cb}\le 2d+1$.
\end{thm}

\end{section}

\begin{section}{Projections onto homogeneous subspaces}

In this section we investigate the complementation of the
homogeneous subspace $\mathcal A_d$ of degree $d$ in the reduced
free product $(\mathcal A, \phi)=*_i(A_i, \phi_i)$. For this
purpose, we put
 $$W_d=\bigoplus_{k=0}^d \mathcal A_k\,,$$
that is, $W_d$ is the closure in $\mathcal A$ of all polynomials
of degree $\le d$. For an element $a$ in the algebraic free
product $A$, we define $\mat P_d(a)$  to be its homogeneous part
of degree $d$. Thus $\mat P_d :   A\to \mat A_d$ is the natural
projection.  Similarly, let  $\mat Q_d : A\mat\to\; W_d$ be the
natural projection.  We also set
$$\mathcal H_d = \bigoplus\limits_{\substack{d\ge n\ge 0 \\
i_1\not=\,i_2\not=\cdots \not=\,i_{n}}} \pl H i n, $$
$$ \mathcal H'_d = \sdd d  \pl H i d.$$

The following is the main result of this section.

\begin{thm}\label{projection}
The natural projection $\mat Q_d$ extends to a completely bounded
projection from $\mathcal A$ onto $W_d$ with $\|\mat Q_d\|_{cb}\le
2d+1$.
\end{thm}

We need some preparations for the proof. Let $\T$ be the unit
circle of the complex plane. For each $z\in \T$ there is a unitary
$U_z$ defined on $\F$ by
 $$\left\{ \begin{array}{ll}
 U_z (\Omega) = \Omega & \\
 U_z (h_1\tens\cdots \tens h_k) = z^k\, h_1\tens\cdots \tens h_k
 & \textrm{ if } k\ge 1
 \textrm{ and } h_1\in \m H_{i_1},\dots ,h_k\in \m H_{i_k}\,.
 \end{array}\right.$$
Note that $U_z^*=U_{\bar z}$. For $n\in \Z$, we let $H_n$ be the
completely contractive projection on $\B(\F)$ defined by
 $$H_n(a) = \int_\T z^n\, U_z^*\, a\, U_z\, dm(z),$$
where the integral is taken with respect to the weak operator
topology, $dm$ being  normalized Lebesgue measure on $\T$. Roughly
speaking, $H_n(a)$ is the part of $a$ that sends tensors of length
$k$ to tensors of length $k+n$. Thus, if $a\in W_d$, then
 $$ a= \sum_{n=-d}^d H_n(a).$$
It would be helpful to understand this formula together with Fact
\ref{factdec}.

\begin{lemma}
Let $a\in \M_s(W_d)$ $(s\in\N)$. Then for $|n|\le d$, 
$$\| (H_n\tens \Id_{\M_s}) (a)\|_{\B(\F)\tens\M_s}\le
\|(H_n\tens \Id_{\M_s}) (a)\|_{\B(\mathcal H_{\big\lceil\frac {d-n}
2\big\rceil},\;\F)\tens\M_s},$$ where for the norm on the right we
consider only restrictions to $\mathcal H_{\big\lceil\frac {d-n}
2\big\rceil}$.
\end{lemma}

\begin{pf}
The point is that for any operator $x\in \B(\F)$,
with respect to the gradation of $\F$, the block matrix of
$H_n(x)=(h_{i,j})_{i,j\ge 0}$ has only non zero coefficients when $i-j=n$,
so
$$\| (H_n\tens \Id_{\M_s}) (a)\|_{\B(\F)\tens\M_s}
 =\sup_{p\ge 1}\|(H_n\tens \Id_{\M_s}) (a)\|_{
 \B(\mathcal H'_p,\;\mathcal H'_{p+n})\tens\M_s}.$$
Thus it suffices to show that this supremum stagnates after
$p=\big\lceil\frac {d-n} 2\big\rceil$. Put $d'=\big\lceil\frac {d-n} 2\big\rceil$ and
let $p>d'$.

Consider the natural partial isometries (obtained by associativity
of tensor products) :
$$ u : \mathcal H'_p \to \mathcal H'_{d'} \tens \F \qquad \textrm{and } \qquad
 v : \mathcal H'_{p+n} \to \mathcal H'_{d'+n} \tens \F.$$
Then, it is easy to check that for any $x\in W_d$
$$ H_n(x)_{|\mathcal H'_p}= v^* (H_n(x)_{|\mathcal H'_{d'}}\tens \Id_{\F}) u.$$
Below is a brief sketch : by linearity, we can
assume that $x$ is homogeneous of degree $k\le d$. By virtue of
the decomposition in Fact \ref{factdec}, we know precisely how $x$
acts as an operator. $H_n(x)$ is the part of $x$ which sends words
of length $p$ to words of length $p+n$. There are two
possibilities :
\begin{itemize}
\item $x$ first annihilates $q$ letters then creates $r$ letters
with $r+q=k$ and $r-q=n$. This is possible only if $k$ and $n$ have the same 
parity and $q=\frac{k-n}2\geq 0$. Here $x$ acts on $q$ letters with $q\leq d'$.

 Of course, we must have $p\ge q$.
\item $x$ first annihilates $q$ letters, then acts once diagonally
and  finally creates $r$ letters with $r+q+1=k$ and $r-q=n$. This is possible only if $k$ and $n$ have different parity and $q=\frac{k-n-1}2$. Here $x$
acts on the first $q+1=\lceil \frac{k-n}2\rceil$ letters of any word.
In
that situation, we must have $p>q$, because there is no diagonal
action on the empty word $\Omega$.
\end{itemize}

In both cases, the maximum number of tensors touched by $x$ is exactly $d'$.
Therefore, acting on a word of length $p$, $H_n(x)$ sees at most
only the first $d'$ letters of the word. This yields the desired
identity from which the lemma  easily follows.

 We also note that from the above discussion, if $d-n$ is odd, $H_n(a)_{|\mat H'_{d'}}$ has a particular shape
as it maps $ \pl H i {d'}$ to $\mathring H_{i_{d'}}\otimes 
\bigoplus\limits_{i_{d'}\not=j_2\not=\,j_2\not=\cdots \not=\,j_{d'+n}}\mathring H_{j_2}\otimes...\otimes\mathring H_{j_{d'+n}}$ (only one diagonal action with $q=d'$).
\end{pf}

\medskip

Fix some $l,k\geq 1$ and let for $i\in I$, $e_i : \mat H'_l \to
\bigoplus\limits_{j_{1}\not=i_2\not=\,j_2\not=\cdots
  \not=\,j_{l-1}\not=i}\mathring H_{j_1}\otimes...\otimes \mathring
H_{j_{l-1}}\otimes \mathring H_i$ be the projection onto words ending
in $\mathring H_i$ and $s_i : \mat H'_k \to
\bigoplus\limits_{i\not=j_{2}\not=i_2\not=\,j_2\not=\cdots \not=\,j_{k}}
\mathring H_i \otimes \mathring H_{j_2}\otimes...\otimes \mathring
H_{j_{k}}$ be the projection onto words starting in $\mathring H_i$.
 For any $x:\mat H'_l\to \mat H'_k$, we let $T(x)=\sum_i s_ixe_i$. Of course
$T$ is complete contraction on $\B( \mat H'_l, \mat H'_k)$.

 We have just noticed that $H_n(a)_{|\mat H'_{d'}}=T( H_n(a)_{|\mat H'_{d'}})$ if $a\in W_d$.

\medskip

\begin{pf}[of Theorem \ref{projection}]

We start by proving  algebraic identities. Let $|n| \le d$,
$d'=\big\lceil\frac {d-n} 2\big\rceil$ and $l\leq d'$.

 We distinguished according to the parity of $d-n$.

If $d-n$ is even for any $a$ in the algebraic free product, we have
$$   H_n(a)_{|\mat H_{d'}} = H_n(\mat Q_d(a)_{|\mat H_{d'}}).$$
This is relevant only if $a \notin W_d$. So assume $a$ is a tensor
of length $k\ge d+1$. We  show that $H_n(a)_{|\mat
H_{d'}}=0$. As above, we have a description of  $H_n(a)$:
\begin{itemize}
\item $a$ first annihilates $q$ letters then creates $r$ letters
with $r+q=k$ and $r-q=n$. Thus $2q=k-n\geq 1+d-n$, and so $q>d'$, this
means that we have to cancel more letters than what we have in
$\mathcal H_{d'}$. \item $a$ first annihilates $q$ letters, then
acts once diagonally and finally creates $r$ letters with
$r+q+1=k$ and $r-q=n$. Thus $2q=k-n-1\geq d-n$, in all case we have
 $q\geq d'$, but then when $x$ would have acted diagonally on
$\Omega$ which is impossible.
\end{itemize}
Therefore, we deduce that $H_n(a)$ vanishes on $\mat H_{d'}$, as
announced.

 If $d-n$ is odd, the situation is slightly more intricate. 
For any $a$ in the algebraic free product, we have 
$$   H_n(a)_{|\mat H_{d'-1}} = H_n(\mat Q_d(a)_{|\mat H_{d'-1}}).$$
And also
$$ T(H_n(a)_{|\mat H'_{d'}}) = H_n(\mat Q_d(a)_{|\mat H'_{d'}}).$$

 The first equality can be treated exactly as above. We focus on the second one.
It is clear if $a\in W_d$ by the observation $H_n(a)_{|\mat H'_{d'}}=T( H_n(a)_{|\mat H'_{d'}})$.  So assume $a$ is a tensor
of length $k\ge d+1$, we  show that $T(H_n(a)_{|\mat
H'_{d'}})=0$. As above, from the description of  $H_n(a)$:
\begin{itemize}
\item $a$ first annihilates $q$ letters then creates $r$ letters
with $r+q=k$ and $r-q=n$. Thus $2q=k-n$. 

If $k>d+1$, then $q\geq 1+\frac{d-n}2$ and $q>d'$ and $a$ cancels more letters than what we have in
$\mathcal H_{d'}$.

If $k=d+1$, then we have a possible non trivial action from $\mathcal H_{d'}$ 
to $\mathcal H_{d'+n}$ correponding to $q=d'$. But then, it will send 
a word ending by $\mathring H_i$ to a word starting by  $\mathring H_j$ with $i\neq j$. this part is killed by $T$.
 \item $a$ first annihilates $q$ letters, then
acts once diagonally and finally creates $r$ letters with
$r+q+1=k$ and $r-q=n$. Thus $2q=k-n-1\geq d-n$, but as $d-n$ is odd
$2q\geq d-n +1$ and  $q\geq d'$, but then when $x$ would have acted diagonally on
$\Omega$ which is impossible.
\end{itemize}

 Thus we have proved the equality.

 In all cases, we can deduce that for $a\in\M_s(A)$ (the algebra of matrices over the algebraic
free product), $|n|\leq d$, 
 $$\|
(H_n\tens \Id_{\M_s}) (Q_d\tens \Id_{\M_s})(a)\|_{\B(\mathcal H_{\big\lceil\frac
{d-n} 2\big\rceil},\;\F)\tens \M_s}\leq \|
(H_n\tens \Id_{\M_s}) (a)\|_{\B(\mathcal H_{\big\lceil\frac
{d-n} 2\big\rceil},\;\F)\tens \M_s}.$$
Indeed this follows immediately from the algebraic identities and 
$$\|
(H_n\tens \Id_{\M_s}) (x)\|_{\B(\mathcal H_{d'},\;\F)\tens \M_s}=\sup_{0\leq l\leq d'} \|
(H_n\tens \Id_{\M_s}) (x)\|_{\B(\mathcal H'_l,\;\F)\tens \M_s}.$$
We conclude using the previous lemma
\begin{eqnarray*}
\| (Q_d\tens \Id_{\M_s}) (a)\|_{\B(\F)\tens \M_s} &\leq & \sum_{n=-d}^d \|
(H_n\tens \Id_{\M_s}) (Q_d\tens \Id_{\M_s})(a)\|_{\B(\F)\tens \M_s} \\&=&
\sum_{n=-d}^d \|
(H_n\tens \Id_{\M_s}) (Q_d\tens \Id_{\M_s})(a)\|_{\B(\mathcal H_{\big\lceil\frac
{d-n} 2\big\rceil},\;\F)\tens \M_s}\\
&\leq &
\sum_{n=-d}^d \|
(H_n\tens \Id_{\M_s}) (a)\|_{\B(\mathcal H_{\big[\frac
{d-n} 2\big]},\;\F)\tens \M_s}\\
 &\leq& \sum_{n=-d}^d \|(H_n\tens \Id_{\M_s}) (a)\|_{\B(\F)\tens \M_s}\\
&\leq & (2d+1) \|a\|_{\B(\F)\tens \M_s}.
\end{eqnarray*}
\end{pf}

\begin{cor}\label{projectionbis}
The natural projection from $\mat P_d: A\to\;\mathcal A_d$ extends
to a completely bounded map on $\mat A$ with norm less than
$\max(4d,\; 1)$.
\end{cor}

This is immediate from Theorem \ref{projection}. In the sequel,
the extensions in Theorem \ref{projection} and Corollary
\ref{projectionbis} will be still denoted by $\mat Q_d$ and $\mat
P_d$, respectively.

\medskip

Corollary \ref{projectionbis} implies that for any $0\le r<1$ the
series
 $$\sum_{k=0}^\infty r^k\mat P_k$$
converges absolutely to a completely bounded map $\mat T_r$ with
cb-norm majorized by
 $$1+4\sum_{k=0}^\infty k\, r^k.$$
However, this estimate for $\|\mat T_r\|_{cb}$ is too bad for
applications. In fact, for the study of approximation properties
in the next section, we will need to know more precisely that
$\mat T_r$ is completely contractive. On the other hand, we will also
need to truncate $\mat T_r$. That $\mat T_r$ is a complete contraction is an
immediate consequence of the following Theorem due to
Blanchard-Dykema \cite{BD}, which will be a main tool for the next
section too.

\begin{thm}\label{dykbl}
Let $(A_i,\phi)$ and $(B_i,\psi_i)$ be $C^*$-algebras with
faithful GNS construction. Let $T_i : A_i\to B_i$ be unital
completely positive maps such that $\psi_i \circ T_i =\phi_i $.
Then there is a completely positive map
$$*_i T_i : *_{i \in I}(A_i,\phi) \to *_{i \in I}(B_i,\psi_i)$$
satisfying, for  $a_1\tens\cdots \tens a_d \in \pl A i d$ with
$\inn i d$,
$$*_i T_i (a_1\tens\cdots \tens a_d) = T_{i_1}(a_1)\tens\cdots
\tens  T_{i_d}(a_d) \in \pl B i d.$$

Moreover, if $A_i$ and $B_i$ are von Neumann algebras and all maps
$T_i$ are normal, then $*_iT_i$ can be extended to a normal map
between the von Neumann free products.
\end{thm}

\begin{prop}\label{truncation}
Let $\mat P_d$ be the natural projection from $\mathcal A$ onto
$\mathcal A_d$ as in Corollary \ref{projectionbis}. Given $0\le
r<1$ and $n\in \N$ define
 $$\mat T_{r}=\sum_{k=0}^\infty r^k \mat P_k
 \quad\mbox{and}\quad
 \mat T_{r,n}=\sum_{k=0}^n r^k \mat P_k.$$
Then $\mat T_{r}$ and $\mat T_{r, n}$ are completely bounded maps
on $\mat A$ with
 $$\|\mat T_r\|_{cb}\le 1\quad\mbox{and}\quad
 \| \mat T_{r,n} \|_{cb}\le 1+ \frac {4nr^n}{(1-r)^2}.$$
The maps $\mat T_{e^{-t}}$, for $t\ge 0$, form a one parameter
semigroup of unital completely positive maps  on $\mat A$
preserving the state $\phi$.

Moreover, the sequence $\mat T_n=\mat T_{(1-1/\sqrt n),n}$ tends
pointwise to the identity of $\mathcal A$ and
 $$\lim_{n\to \infty}\|\mat T_n\|_{cb}=1.$$
\end{prop}

\begin{pf}
 We apply  Theorem \ref{dykbl} with
$(B_i,\psi_i)=(A_i,\phi_i)$
 and $T_i(a)=U_{r,i}(a)=r a + (1-r)\phi_i(a) 1$. These maps are obviously
unital, completely positive and  preserve the states. Formally,
 $$*_iU_{r,i} = \sum_{k=0}^\infty r^k \mat P_k.$$
There does not exist any trouble since the series above is
absolutely convergent for any $0\le r<1$, as already observed
previously. Moreover, we know that $\|*_i U_{r,i}\|_{cb}= 1$. Thus
by the triangular inequality:
 $$\| T_{r,n} \|_{cb}\le \| *_i U_{r,i} \|_{cb}+\sum_{k=n+1}^\infty
 r^k\|\mat P_k\|_{cb} \le  1+ \frac {4nr^n}{(1-r)^2}.$$
As $*_iU_{r,i}$ is bounded uniformly in $r$ and tends to the
identity pointwise on $\mathcal A$ as $r\to 1$, the second
assertion follows from simple computations as $\mat T_n$ is just a
perturbation of $*_iU_{(1-1/\sqrt n),i}$.
\end{pf}

\medskip

The semigroup $(\mat T_{e^{-t}})_{t\ge0}$ above is called the
Poisson semigroup or kernel on $\mat A$ because of its analogy
with the usual Poisson kernel on the unit circle.

\medskip

\begin{rk}\label{normal projection}
All previous results extend to the case of amalgamated reduced
free products or von Neumann algebra reduced free products. In the
latter case, all maps constructed above are normal. Just note that
the normality of $\mat P_d$ follows from that of the Poisson
kernel $\mat T_r$.
\end{rk}

\medskip

\begin{rk}
The results analogous to those presented in this section are well
known for free groups and are due to Haagerup \cite{Haa}.
Corollary \ref{projectionbis} in the case of free groups can be
also deduced from the fact that a free group acts non trivially on
a tree \cite{Boz}. This is also related to another property;
the Hertz-Schur multiplier $z^{|g|}$ on a free group is a
coefficient of a uniformly bounded representation. This latter
result is due to Bo{\.z}ejko \cite{Boz2}. If the $A_i$'s are reduced group
$C^*$-algebras (say $A_i=C^*_\lambda(G_i)$),
 the projection $\mat P_d$ is an Hertz-Schur multiplier on
the group free product $*_{i=1}^NG_i$. Following a previous
unpublished work of Haagerup and Szwarc,  Wysocza{\'n}ski
\cite{Wy} was able to compute explicitly the norm of such
multipliers in terms of trace norms of some Hankelian matrices.
Using his result, one can easily get that if $N=\infty$ and all
the groups are infinite then $\|\mat P_d\|_{cb}/d\to_{d\to \infty}
8/\pi$. This shows that the linear growth estimate is the best
possible.

\end{rk}

\end{section}

\begin{section}{Applications}

We will apply the results in the previous two sections to study
various approximation properties for reduced free products.

\subsection{Exactness}

 An operator space $X$ is said to be exact, if there is a
constant $\lambda$ such  that  for any (closed) ideal $I$ of a
$C^*$-algebra $B$, the natural map
$$T :X \tens_{\min} B/(X \tens_{\min} I) \to X \tens_{\min} B/I$$
is invertible with $\|T^{-1}\|\le \lambda$. The exactness constant
of $X$ is then the infimum over all such constants $\lambda$.

There are various equivalent reformulations of exactness, we refer
to \cite{Pis} for more information. Among well- know results, we
mention that if a $C^*$-algebra $A$ is exact, then its exactness
constant is 1. Moreover, the $\lambda$-exactness property can be
thought as an approximation property : the inclusion of $A$ into
$\B(H)$ can be approximated in the point norm topology by finite
rank unital completely positive maps, with factorization norm
through $\M_n$ bounded by $\lambda$.

The column and row Hilbert spaces, the space $\K$ of compact
operators on $\ell_2$ are examples of 1-exact operator spaces.
More generally, $X$ and $\K\tens_{\min} X= C\tens_h X\ \tens_h R$
have the same exactness constant, where $C$ and $R$ are
respectively the row and column space based on $\ell_2$. We will
also need that a direct sum (in the $\ell_\infty$-sense) of
$\lambda$-exact operator spaces is also $\lambda$-exact.

Stability of exactness under (amalgamated) reduced free product
was first proved by Dykema \cite{D} and another proof was later
given  by Dykema-Shlyakenthko \cite{DS}.

Here we present a simpler proof using our Khintchine inequalities.
It consists of a mere adaptation of a nice argument by Pisier (see
chapter 17 in \cite{Pis}).

\medskip

Let $(A_i,\phi_i)_{i\in I}$ be $C^*$-algebras with faithful GNS
constructions and ($\mathcal A, \phi)=*_{i\in I} (A_i,\phi_i)$.

\begin{thm}
 If all $C^*$-algebras $A_i$ are exact, then $\mathcal A$ is exact.
\end{thm}

\begin{pf}
Since we are dealing with $C^*$-algebras, we have to show that for
any ideal $I\lhd B$ of a $C^*$ algebra $B$ and any element $x$ in
$\mathcal A\tens B$, we have
$$\| (\Id\tens \rho) x\|_{\mathcal A \tens_{\min} B/I}=
 \| \tilde \rho  (x)\|_{(\mathcal A \tens_{\min} B)/(\mathcal A \tens_{\min}
 I)},$$
where $\rho:B\to B/I $ and $\tilde\rho : \mathcal A \tens_{\min} B\to
(\mathcal A \tens_{\min} B)/(\mathcal A \tens_{\min} I)$ are the quotient maps.

According to Theorem \ref{khintchine}, all the subspaces $\mathcal
A_d$ are exact with constant less than $2d+1$ for $E_d$ is exact
with constant 1. Since $\mathcal A_d$ is $4d$-complemented in
$\mathcal A$ by virtue of Corollary \ref{projectionbis}, we deduce
that $W_d=\oplus_{k\le d} \mathcal A_k$ (algebraically) is also
exact. We just need to compare $W_d$ with $\ell_\infty ((\mathcal
A_k)_{k\le d})$, so we obtain an exactness constant bounded by
$4d\,(\sum_{k\le d} 2k+1)\leq 4(d+1)^3$. In turn, this means that
for any $x\in W_d\tens B$, we have
$$4(d+1)^3\| (\Id\tens \rho) x\|_{\mathcal A \tens_{\min} B/I}\geq
 \| \tilde \rho  (x)\|_{(\mathcal A \tens_{\min} B)/(\mathcal A \tens_{\min}
 I)}.$$
Now, we can use the usual trick that $\|x\|^{2n}=\|(x^*x)^n\|$. If
we start with $x\in W_d\tens B$, then  $(x^*x)^n$ belongs to
$W_{2nd}\tens B$. As $\rho$ and $\tilde \rho$ are
$*$-representations, we get by the previous observations :
$$ 4(2dn+1)^3\| (\Id\tens \rho) x\|^{2n}_{\mathcal A \tens_{\min} B/I}\ge
 \| \tilde \rho  (x)\|^{2n}_{(\mathcal A \tens_{\min} B)/(\mathcal A
\tens_{\min} I)}.$$ Taking $2n$-root and letting $n\to \infty$, we
obtain the desired inequality for any $x\in W_d\tens B$ and every
$d$. The conclusion follows by the density of  the $W_d$'s in
$\mathcal A$.
\end{pf}

\medskip

\begin{rk}
 Since Khintchine inequalities are true for amalgamated free product, the same
proof can be carried out to the amalgamated case (see the last
section for more details).
\end{rk}

\subsection{CCAP}

An operator space $X$ is said to have the completely bounded
approximation property (CBAP) with constant $\lambda$ if the
identity of $X$ is the limit for the pointnorm topology of a net
of finite rank $\lambda$-completely bounded maps. As far as we are
concerned, we will deal with the completely contractive
approximation property (CCAP), that is the CBAP with constant 1.

Similarly, if $X$ is a dual space, it has the weak-$*$ CCAP if the
identity of $X$ is  the limit in  the point weak-$*$ topology of a
net of finite rank completely contractive weak-$*$ continuous
maps. This is equivalent to say that the predual $X_*$ of $X$  has
the CCAP.

For $C^{*}$-algebras, the CCAP is a notion stronger than
exactness. It is an open question to know if the CCAP passes to
reduced free products. In this section, we give some  partial
answers.

\medskip

Our first result on the CCAP is the following

\begin{thm}
If  all the $A_i$ are finite dimensional, then their reduced free
product $\mathcal A$ has the CCAP.
\end{thm}

\begin{pf}
This is an immediate consequence of Proposition \ref{truncation}
for  the maps $\mat T_n$ defined there are of finite rank in the
present situation.
\end{pf}

\medskip

\begin{rk}\label{state preserving to ccap}
More generally, assume that for each $i$, there is a net of finite
rank unital completely positive maps preserving the state $\phi_i$
on $A_i$ converging to the identity pointwise in norm, then
$\mathcal A$ has the CCAP. This follows directly by composing the
truncated Poisson kernel with the free product of the
approximating sequences. This condition implies, of course,  that
the $A_i$ are nuclear.
\end{rk}

Thus the problem of proving the CCAP for a reduced free product
reduces to find an approximating net as in the remark above for
each factor $C^*$-algebra in the product. With this in mind we
have the following

\begin{prop}\label{state preserving}
Let $A$ be a nuclear unital $C^*$-algebra and $\phi$ a faithful
state on $A^{**}$. Then
there is a net of finite rank unital completely positive maps
$T_i:A\to A$ converging to the identity pointwise in norm such
that $\phi\circ T_i=\phi$.
\end{prop}

\begin{pf}
Since $A$ is nuclear, $A^{**}$ is semidiscrete. So there are
normal unital completely positive maps $V_i: A^{**}\to A^{**}$
converging pointwise to the identity for the weak-$*$ topology,
admitting  a normal completely positive factorisation $V_i=\beta_i
\alpha_i$ with
 $$\alpha_i : A^{**}\to \M_{d_i}, \qquad \beta_i : \M_{d_i}\to  A^{**}
 \qquad
 \textrm{ and }\qquad  \lim \|\alpha_i\|=\lim \|\beta_i\|=1.$$

It is well known that completely positive maps $\M_d\to B$ with
values in a $C^*$-algebra $B$ are identified with positive
matrices with values in $B$. By Kaplansky's density theorem, we
can assume (by passing to another net) that the maps $\beta_i$ are
actually with values in $A$, and moreover such that the $V_i$ are
still unital. Let us justify the last point. Fix $i$. Let
$(\alpha_i(j,k))\in \M_{d_i}$ be the matrix of $\alpha_i(1)$ and
$(\beta_i(j,k))\in \M_{d_i}(A^{**})$ the matrix corresponding to
$\beta_i$. Since $V_i$ is unital,
 $$\sum_{j,k=1}^{d_i}\alpha_i(j,k)\,\beta_i(j,k)=1.$$
Let $\beta_i^\nu=(\beta_i^\nu(j,k))\in \M_{d_i}(A)_+$ be a net
converging to $\beta_i$ in the weak-$*$ topology. Then
 $$a^\nu\,\mathop{=}^{\textrm{def}}\,
 \sum_{j,k=1}^{d_i}\alpha_i(j,k)\,\beta_i^\nu(j,k)\to 1
 \quad\textrm {weakly in}\ A.$$
Thus passing to convex combinations, we can assume $a^\nu\to 1$ in
norm. Consequently, $\|a^\nu\|\to 1$. Then dividing by
$\|a^\nu\|$, we can further assume $0\le a^\nu\le 1$. Now we
define $\widetilde \beta_i^\nu: \M_{d_i}\to  A^{**}$ by
 $$\widetilde \beta_i^\nu(m)
 =\sum_{j,k=1}^{d_i}m(j,k)\,\beta_i^\nu(j,k)
 +\frac{\sum_{j=1}^{d_i}m(j,j)}{\sum_{j=1}^{d_i}\alpha_i(j,j)}\,
 (1-a^\nu),\quad m\in \M_{d_i}.$$
Then we still have that $\widetilde \beta_i^\nu\to \beta_i$ in the
weak-$*$ topology; but now $\beta_i^\nu\circ\alpha_i$ is unital.

Let  $V_{i*}$ be the pre-adjoint of $V_i$. Then
$V_{i*}=\alpha_{i*}\beta_{i*}$. Taking convex combinations, the
maps $V_{i*}$ can be assumed to tend to the identity pointwise in
norm on $A^*$.

 At this stage, we have that $\|\phi\circ V_i-\phi\|\to 0$. By the
Jordan decomposition there are positive linear functionals
$\psi_i, \xi_i$ in $A^*$ such that $\phi\circ V_i-\phi =
\psi_i-\xi_i$ and $\|\psi_i\|+\|\xi_i\|\to 0$. Let
$U_i(x)=V_i(x)+\xi_i(x) 1$. It is easy to check that the sequence
$(U_i)$ shares the  properties of $(V_i)$ : it consists of  normal
 completely positive maps tending to the identity
pointwise in the weak-$*$ topology such that $\lim \|U_i\|=1$ and
$U_{n*}$ tends pointwise in norm to the identity on $A^*$. The
point now is that as functionals
 $$\phi\circ U_i \ge \phi.$$
Hence by the Radon-Nikodym theorem, there are selfadjoint elements
$a_i$ in $A^{**}$ such that
 $$0\le a_i \le 1\qquad \textrm{and} \qquad \forall \, x\in A^{**},
 \;\;\phi(U_i(a_ixa_i))=\phi(x).$$
Passing to a subnet if necessary, we can assume that $a_i$ and
$a_i^2$ converge to some $a$ and $b\in A^{**}$ $*$-weakly.
Clearly, $0\le b\le a\le 1$. First, we have
$\phi(U_i(a_i^2))=\phi(1)=1$. On the other hand, as $U_{n*}$ tends
to the identity pointwise in norm and the $a_i$ are uniformly
bounded, we conclude that $U_i(a_i^2)$ tends to $b$ for the
weak-$*$ topology. So $\phi(b)=1$, and $a=b=1$ by the faithfulness
of $\phi$. Moreover, this shows that actually $a_i\to 1$ for the
strong operator topology of $A^{**}$ (in its universal
representation).

Now we define $S_i: A^{**}\to A$ by
 $$S_i(x)=U_i(a_ixa_i)$$
Then $S_i$ is a finite rank completely positive maps preserving
the state $\phi$ on $A^{**}$. Moreover, the same kind of arguments as above
leads to the fact that the net $S_i$ converges to the
identity of $A^{**}$ weak-$*$-pointwise with $\lim \|S_i\|=1$.
Restricting to $A$ and taking some convex combinations, we get
finite rank completely positive maps $S'_i: A\to A$ converging to
the identity in norm and preserving the state with $\lim
\|S'_i\|=1$. Finally, the desired net $T_i$ is given by
 $$T_i(x)= \frac 1 {\|S'_i(1)\|} S'_i(x)
 +\phi(x)\Big[1- \frac 1 {\|S'_i(1)\|}S'_i(1)\Big].$$
We have thus achieved the proof.
 \end{pf}

\medskip

In the case of von Neumann algebras, one can easily adapt the
arguments in the proof of Proposition \ref{state preserving} to
get the

\begin{prop}\label{vnstate preserving}
Let $M$ be a semidiscrete von Neumann algebra and $\phi$ a normal
state on $M$. Then there are normal finite rank unital completely
positive maps preserving $\phi$ and converging to the identity in
the point weak-$*$ topology.
\end{prop}

\begin{pf}
If $\phi$ is faithful, then the result follows from an easy
adaptation of the proof of Proposition  \ref{state preserving}.
 \\
Let us consider the general case. \\
 Let $p$ be the support
 of $\phi$. Assume that the central support $c(p)$ of $p$ is 1.
 Then by the comparison theorem for projections
(see Lemma 1 in Chapter III paragraph 2 of \cite{Dix}) and Zorn's
lemma, there is a familly of orthogonal projections $(p_i)_{i \in
I}$ such that
 $$
 \textrm{i)}\  p=p_{i_0} \textrm{ for some } i_0 \in I; \quad
 \textrm{ii)}\  p_i  \preccurlyeq p;\quad
 \textrm{iii)}\ \sum_{i\in I} p_i =1.
 $$
Then it follows that for $H=\ell_2^I$ with an orthonormal basis
$(e_i)_{i\in I}$, there is a projection $q\in pMp\overline
\tens\,\B(H)$ with $q\geq p\tens P_{e_{i_0}}$, so that
 $$ M \simeq q (pMp \,\overline\tens \, \B(H))q.$$
With this identification the state $\phi$ is the restriction of
$\phi\tens \omega_{e_{i_0}}$ to $M$.

 Since $pMp$ is semidiscrete and $\phi$ is faithful on it,
we already know that there is a net $U_\alpha: pMp\to pMp$ of
normal finite rank unital completely positive maps preserving
$\phi$ that converges to the identity for the weak-$*$ topology.
Similarly, it is not hard to construct a family of normal finite
rank unital completely positive maps $V_\alpha$ on $\B(H)$
preserving $\omega_{e_{i_0}}$ that converges to the identity for
the weak-$*$ topology on $\B(H)$.
 Then, the maps $T_\alpha: M\to M$ given by
  $$T_\alpha(x)=q(U_\alpha\overline \tens\,V_\alpha)(x)q$$
have the required properties.

Finally, if $c(p)\neq 1$, we have $M= c(p) M \oplus c(p)^\bot
M$. \\ To construct the approximating net, we use the above argument
for the $ c(p) M$ part and the semidiscretness of  $c(p)^\bot M$
(since $\phi$ is zero on it).
\end{pf}

\medskip

\begin{rk} The argument for the reduction of the general case to the
faithful case was kindly shown to us by Uffe Haagerup.
\end{rk}

\medskip

The preceding proposition implies a corresponding result on the
weak-$*$ CCAP for von Neumann reduced product. To this end let us
first recall the following well-known elementary fact. Let $B$ be
a $C^*$-algebra with a given state $\rho$ and $V$ a unital
completely positive map on $C$ such that $\rho\circ V=\rho$. Then
$V$ naturally induces a contraction on $L_2(B,\rho)$ (resp.
$L_2(B,\rho)^{op}$). The same is true if $B$ is a von Neumann
algebra and $\rho$, $V$ are normal.

\begin{thm}
 Let $(M_i,\phi_i)$ be semidiscrete von Neumann algebras with
distinguished normal states with faithful GNS constructions. Then
the reduced von Neumann free product $\overline *_i (M_i,\phi_i)$
has the weak-$*$ CCAP.
\end{thm}

\begin{pf}
We keep the notations at the end of section~2. Once again, we use
Theorem \ref{dykbl} and the von Neumann version of Proposition
\ref{truncation} to get a normal Poisson kernel on $\overline{\mat
M}$, where $\overline{\mat M}= \overline *_i (M_i,\phi_i)$. Note
that all the projections $\mat P_d$ (onto the weak-$*$ span of
$\Sigma_d$) are normal (see Remark \ref{normal projection}). By
Proposition \ref{vnstate preserving}, for each $i$, we get a net
of finite rank unital state preserving completely positive maps
$(U_{i,n})_n$ converging to the identity for the point-weak-$*$
topology on $M_i$. Consequently, these maps give rise  to finite
rank contractions on $L_1^m$ and  $K_1^{k}$ converging to the
identity for the weak topology (as Hilbert spaces). Since $L_1^m$
and  $K_1^{k}$ are homogeneous operator spaces, these contractions
are automatically completely contractive. In particular, for each
$d$, on $\overline X_d$, their natural extensions (tensorization
on each component) $G_{d,n}$  tends $*$-weakly to the identity.

We consider approximating sequences obtained by the previous
scheme, the truncations $\mat T_N$ of the Poisson kernel composed
with the free product maps $(U_n=*U_{i,n})_n$. We obtain a net of
normal unital finite rank completely bounded maps $V_{N,n}=U_n
\mat T_N : \overline{\mat M}\to\;  \overline{\mat M}$. It remains
to verify that it converges $*$-weakly to the identity. For fixed
$N$, $V_{N,n}$ factorizes normally through $\oplus_{k\le N}
\overline X_d$ by virtue of  the Khintchine inequalities. We have
a commutative diagram (with $r=1-1/\sqrt N$):
$$\xymatrix{
\overline{\mat M}\ar[rr]^{V_{N,n}} \ar[d]_{\sum_{d\le N} r^d
\Theta_d \mat P_d}
& &\overline{\mat M}\\
\oplus_{d\le N} \overline X_d\ar[rr]^{\oplus_{k\le N} G_{d,n}}& &
\oplus_{d\le N} \overline X_d \ar[u]_{\sum_{d\le N} \Xi_d}}$$
 As $G_{d,n}$ tends $*$-weakly to the identity and all maps in the
diagram are normal, we get that for any $x\in \overline{\mat M}$
$$\lim_n U_n \mat T_N(x)= \mat T_N(x).$$
 So, by general principles, it follows that
$$\lim_N\lim_n U_n \mat T_N(x)=x.$$
\end{pf}

\medskip

At the time of this writing we do not know whether any nuclear
unital $C^*$-algbera $A$ possesses a net of finite rank completely
positive unital maps preserving a given state. It is however easy
to see that the answer is affirmative for the unitization $\tilde
\K$ of the compact operator algebra $\K$. Indeed, given a state
$\phi$ on $\K$, considering a basis that diagonalizes $\phi$ and
using truncation with respect to that basis gives a net of maps on
$\K$ with all requirements . (Note that on $\K$, the GNS
representation of a state is always faithful.)
\begin{prop}
For any states $\psi, \phi$ on $\tilde \K$, the reduced
free product $(\tilde \K,\phi)*(\tilde \K,\psi)$ has the CCAP.
\end{prop}

The nuclearity/semidiscretness hypothesis actually fits very well
with our situation. Next, we try to drop this hypothesis. The very
first obstacle we encounter is how to extend a completely bounded
state preserving map on a $C^*$-algebra equipped with a state to
the associated $L_2$-space. In fact, if we want to approximate the
identity on a reduced free product, the operator spaces appearing
in the Khintchine inequalities suggest that we should perform
simultaneously approximations for each algebra  and the associated
$L_2$ spaces. If we assume only the CCAP for the $C^*$-algebra, it
is not even clear that the approximating sequence can be chosen to
be bounded in the $L_2$-norm with respect to any state. In the
following lemma, we will assume that all endomorphisms in
consideration have natural bounded extensions to the corresponding
$L_2$-spaces.

\begin{lemma}\label{prod}
Let $d>0$ be fixed and for any $1\le k \le d$ and $i\in I$, let
$T_{i,k}: A_i\to A_i$ be  completely bounded maps such that
$\phi_i\circ T_{i,k}=\phi_i$ and $T_{i,k}$ naturally extends to
bounded maps on $L_2(A_i,\phi)$ and $L_2(A_i,\phi)^{op}$. Then the
map $\prod_k T_{i,k}:\Sigma_d \to \Sigma_d$ defined by,   for
$a_1\tens\cdots \tens a_d \in \pl A i d$ with $\inn i d$,
$$(\prod_k T_{i,k}) (a_1\tens\cdots \tens a_d) =
T_{i_1,1}(a_1)\tens\cdots \tens T_{i_d,d}(a_d) \in \pl A i d,$$
admits a completely bounded extension on $\mathcal A_d$ with
cb-norm majorized by
 $$(2d+1) \prod_{k=1}^d
\max_i \big\{ \max \{\|T_{i,k}\|_{cb(A_i,A_i)},
\|T_{i,k}\|_{\B(L_2(A_i),\, L_2(A_i))},
\|T_{i,k}\|_{\B(L_2(A_i)^{op},\, L_2(A_i)^{op})}\}\big\}.$$
\end{lemma}

\begin{pf}
 This map $\prod_k T_{i,k}$ is well defined algebraically,
as $\m A_i$ is stable by $T_{i,k}$ for any $k$. For each $k$, the
direct sum $\oplus_i\, T_{i,k}$ extends to completely bounded
maps, say, $G_k$ on $\ell_\infty ((A_i))$, $V_k$ on $L_1$ and
$U_k$ on $K_1$. On each component of $X_d$, we can consider the
tensor product map of them. We thus obtain a completely bounded
endomorphism of $X_d$, which, somehow, is an extension of $\prod_k
T_{i,k}$ defined on $\iota(\Sigma_d)$. The conclusion is then
obvious using Theorem \ref{khintchine} and norm estimates for
tensor products.
\end{pf}

\medskip

Now assume that $( B_i,\psi_i)$ are unital $C^*$-algebras with
distinguished states $\psi_i$ with faithful GNS constructions. Let
$A_i\subset B_i$ be unital $C^*$-subalgebras such that all
$\phi_i={\psi_i}_{|A_i}$ also have faithful GNS constructions.

Suppose that for each $i$, there is a net of finite rank maps
$(V_{i,j})_{j\in J}$ on $A_i$ converging to the identity
pointwise, preserving the state and such that $\limsup_j
\|V_{i,j}\|_{cb}=1$. (consequently, all $A_i$ have the CCAP).

 Assume moreover that for each  pair $(i,j)$, there is a completely positive unital
map $U_{i,j}: A_i \to B_i $ preserving the states, such that
$$\lim_j \|V_{i,j}-U_{i,j}\|_{cb}+
\|V_{i,j}-U_{i,j}\|_{\B(L_2(A_i,\phi_i),L_2(B_i,\psi_i))}+
\|V_{i,j}-U_{i,j}\|_{\B(L_2(A_i,\phi_i)^{op},L_2(B_i,\psi_i)^{op})}=0.$$

\begin{prop}
 Under the above hypothesis, $\mathcal A$ has the CCAP.
\end{prop}

\begin{pf}
For simplicity we can assume that $I$ is finite, which is not
relevant for the CCAP.

The first assumption allows us to view $\mathcal A$ as a
subalgebra of $\mathcal B$ by \cite[Proposition 2]{BD}.

Since the maps $U_{i,j}$ are completely positive and preserve the
states, they extend to contractions from $L_2(B_i,\psi_i)$ to
$L_2(A_i,\phi_i)$ and from $L_2(B_i,\psi_i)^{op}$ to
$L_2(A_i,\phi_i)^{op}$.  Passing to a subnet if necessary, we may
assume that $\epsilon_{i,j}\le 1$ for all $i$ and $j$, where
 $$\epsilon_{i,j}=\|V_{i,j}-U_{i,j}\|_{cb} +
 \|V_{i,j}-U_{i,j}\|_{\B(L_2(A_i,\phi_i),L_2(B_i,\psi_i))}+
 \|V_{i,j}-U_{i,j}\|_{\B(L_2(A_i,\phi_i)^{op},
 L_2(B_i,\psi_i)^{op})}.$$
Then we deduce that the $V_{i,j}$ also extend to bounded maps on
$L_2(A_i,\phi_i)$ and $L_2(A_i,\phi_i)^{op}$ and
 $$\sup\big\{\|V_{i,j}\|_{cb}\,,\;
 \|V_{i,j}\|_{\B(L_2(A_i,\phi_i))}\,,\;
 \|V_{i,j}\|_{\B(L_2(A_i,\phi_i)^{op})},\ i\in I, j\in J\big\}\le 2.$$

Thus Lemma \ref{prod} yields the product maps $F_{d,j}=\prod_k
V_{i,j} : \mathcal A_d \to\; \mathcal A_d$ (without any dependence
on $k$ here). To get a finite rank approximation of the identity
of $\mathcal A$, we consider the following (with $r=1-1/\sqrt N$)
 $$ D_{N,j}= \sum_{d=0}^N r^d F_{d,j} \mat P_d.$$
 These maps are all of finite rank and completely bounded. Obviously,
$(D_{N,j})$ tends (as $N,j\to \infty$) to the identity on the
dense vector subspace of $\mathcal A$ consisting of linear
combinations of elementary tensors. To conclude, we only need to
check that $\lim_{N} \lim_j \|D_{N,j}\|_{cb}=1$.

To that end, by Theorem \ref{dykbl}, we get unital completely
positive maps $*_i\,U_{i,j} :\mathcal A\to \mathcal B$. In
$\mathcal B$, we can consider the approximation of the identity
$\mat T_N$ given by Proposition \ref{truncation}. Define  the maps
 $$E_{N,j}=*_i\,U_{i,j}\circ \mat T_N=\mat T_N\circ
 *_i\,U_{i,j}.$$
Then we have $\limsup_{j,N} \|E_{N,j}\|_{cb}=1$. By the triangular
inequality
 $$\|D_{N,j}\|_{cb} \le \|E_{N,j}\|_{cb} +\|E_{N,j}-D_{N,j}\|_{cb}.$$
We have
 $$\|E_{N,j}-D_{N,j}\|_{cb}\le \sum_{d=0}^N \| (F_{d,j}-
 *_i\,U_{i,j})\mat P_d\|_{cb}.$$
However, $*_i\,U_{i,j}=\prod_k U_{i,j}$ on $\mathcal A_d$ as in
Lemma \ref{prod}. Thus we deduce:
\begin{eqnarray*}
 \|E_{N,j}-D_{N,j}\|_{cb}
 &\le&\sum_{d=1}^N 4d\, \big\|\prod_k V_{i,j}-\prod_k U_{i,j}\big\|_{cb}\\
 &\le&\sum_{d=1}^N 4d\,\big\|\sum_{l=1}^{d}\prod_{k=1}^{l-1}V_{i,j}
 \prod_{k=l}^d U_{i,j}- \prod_{k=1}^{l}V_{i,j}
 \prod_{k=l+1}^d U_{i,j}\big\|_{cb}\\
 &\le& \sum_{d=1}^N 4d\, \sum_{l=1}^{d}\big\| \prod_{k=1}^{l-1} V_{i,j}
 (V_{i,j}-U_{i,j})\prod_{k=l+1}^d U_{i,j}\big\|_{cb}\\
 &\le& \sum_{d=1}^N 4d\,\sum_{l=1}^{d} 2^{l-1}\max_{i}
 \epsilon_{i,j}\,.
\end{eqnarray*}
Hence
$$\lim_N \lim_j \|  D_{N,j} \| \le 1.$$
\end{pf}

\medskip

For group $C^*$-algebras, we have more information about the CBAP.
We refer to \cite{Boz,BP,KH} for background on the subject. Given
a discrete group $G$, we denote as usual by $C^*_\lambda (G)$ and
$VN(G)$ the reduced $C^*$-algebra and von Neumann algebra of $G$,
respectively. These algebras are generated by the left regular
representation $\lambda$ of $G$ on $\ell_2(G)$ and equipped with
the canonical tracial faithful vector state $\omega_{\delta_e}$.
Here, we use the standard notation $\{\delta_g \,;\, g\in G\}$ for
the canonical basis of $\ell_2(G)$ ($e$ being the neutral element
of $G$).

If $f:G\to \C$, we denote by $M_f: C^*_\lambda (G) \to C^*_\lambda
(G)$ the Hertz-Schur multiplier associated with $f$. This is the
linear map defined by $M_f(\lambda(g))= f(g)\lambda(g)$ for $g\in
G$. If $a=(a_{s,t})_{s,t\in G}$ is a matrix index by $G$,  $M_a$
denotes the Schur multiplier on $\B(\ell_2 (G))$ induced by $a$ in
the natural basis. So, we have the relation $M_f=M_a$ with
$a=(f(s^{-1}t))$.

A discrete group $G$ is weakly amenable if there is a net of
functions $f_i:G\to \C$  with finite support such that $f_i$
converges pointwise to the constant function 1 and such that
$\limsup_i \|M_{f_i}\|_{cb}<\infty$. The Haagerup constant
$\Lambda(G)$ of $G$ is then defined to be the infimum over all
such $\limsup$.

Let us recall the following well-known result \cite{Boz,KH}.

\begin{thm}
Let $G$ be a discrete group and $C>0$. Then the following
properties are equivalent
\begin{enumerate}[\rm i)]
 \item $G$ is weakly amenable with constant $C$.
 \item $C_\lambda^*(G)$ has the CBAP with constant $C$.
 \item $VN(G)$ has the weak-$*$ CBAP with constant $C$.
\end{enumerate}
\end{thm}

Haagerup first proved that the free groups are examples of groups
with $\Lambda(G)=1$ even if they are not amenable. This result was
generalized by Bo{\.z}ejko and Picardello in \cite{BP} to the free
product of amenable groups (with amalgamation over a finite
group). This can be deduced quite easily from our previous
results. The connection between free product and group algebras is
quite simple for the reduced $C^*$-algebra of the free product of
groups is nothing but the reduced free product of the group
$C^*$-algebras with respect to their usual tracial states.

\begin{thm}
Let $(G_i)_{i\in I}$ be weakly amenable discrete groups with
constant 1. Then $\mathop{*}\limits_{i\in I} G_i$ is also weakly
amenable with constant 1.
\end{thm}

\begin{pf}
We just need to check the assumption of the previous proposition
for a single weakly amenable group, with $A=C_\lambda^*(G)$,
$B=\B(\ell_2(G))$ and $\psi=\omega_{\delta_e}$ (the GNS
construction is then the identity).

Since $G$ is weakly amenable with constant 1, there is a net of
functions $f_i:G\to \C$ with finite support converging pointwise
to the constant $1$ and such that $\lim_i \|M_{\tilde
f_i}\|_{cb}=1$, where $ M_{\tilde f_i}$ is the Schur multiplier on
$\B(\ell_2(G))$ with symbol $\tilde f_i$ given by $\tilde
f_i=(f_i(s^{-1}t))_{s,t\in G}$.

Without loss of generality, we can assume that $f_i$ is real
(replacing it by $(f_i+\overline f_i)/2$), symmetric (i.e.
$f_i(s)=f_i(s^{-1})$ for any $s\in G$) and $f_i(1)=1$. Thus,
$M_{\tilde f_i}$ is represented by a selfadjoint real matrix and
preserves the state.

Suppose that $\|M_{f}\|\le 1+ \epsilon$ for some $\epsilon>0$ (we
drop the index $i$). From the representation of Schur multipliers
(see, e.g. \cite{Pis}), we know that there exist two families
$(x_g)_{g\in G}$ and  $(y_g)_{g\in G}$ of vectors in some Hilbert
space $H$, such that
 $$\sup_g \|x_g\|= \sup_g \|y_g\| \le \sqrt{1+\epsilon} \qquad
 \textrm{and} \qquad
 \forall s,t\in G, \quad f(s^{-1}t)=\langle x_s,\;y_t\rangle.$$
Consider the Schur multipliers with symbols
 $$a=\big(\langle\frac {x_s+y_s}2,\; \frac {x_t+y_t} 2\rangle\big)_{s,t}
 \qquad  \textrm{ and } \qquad
 b=\big(\langle\frac {x_s-y_s} 2,\; \frac {x_t-y_t} 2\rangle\big)_{s,t}.$$
Then $M_a$ and $M_b$ are completely positive Schur multipliers and
the assumption on $f$ gives (by polarization) that $M_{\tilde
f}=M_a-M_b$.
  However,
 $$\|M_a\|_{cb}=\|M_a(1)\|=\sup_s \left\|\frac {x_s+y_s} 2\right\|^2\le
 1+\epsilon,$$
so computing on 1 yields that for any $s\in G$
  $$1 + \left\|\frac {x_s-y_s} 2\right\|^2=\left\|\frac {x_s+y_s} 2
 \right\|^2\le 1+\epsilon;$$
whence  $\|M_b\|_{cb}= \sup_s \|\frac {x_s-y_s} 2\|^2\le\epsilon$.
Consequently, $M_{\tilde f}$ is an $\epsilon$-perturbation of the
completely positive map $M_a$. To conclude, it suffices to notice
that there is a  positive diagonal matrix $d$ such that $\frac 1
{1+\epsilon}M_{a+d}$ is unital completely positive (so preserves
the state for it is a Schur multiplier) with $\|M_d\|_{cb}\le
\epsilon$. Thus we get the desired approximation at least for the
cb-norm. The part for the $L_2$ estimates are also satisfied
because any contractive Schur multiplier defines a contraction on
$L_2(\B(\ell_2(G)),\omega_{\delta_e})=\ell_2(G)$ (and similarly
for the opposite).
\end{pf}

\medskip

For group von Neumann algebras and more generally finite von
Neumann algebras, a yet more approximation property has attracted
a lot of attention since many years. It is the so-called Haagerup
property: A finite (separable) von Neumann algebra $M$ with a
faithful normal trace $\tau$ has the Haagerup property if there is
a sequence of unital completely positive maps $T_n:M\to M$ such
that $\tau\circ T_n=\tau$, whose natural extensions to $L_2$ are
compact and $(T_n)$ converges to the identity for the point norm
topolgy in $L_2$. Of course, this implies that $T_n$ actually
tends to the identity of $M$ for the point weak-$*$ topology.

It is know that this property is preserved by reduced free product
(\cite{Boca} and \cite{J}). The main difference with the CCAP
relies on the fact that the map $T_n$ are only compact and not of
finite rank. It is however conjectured that a group is weakly
amenable with constant 1 if and only if it has the Haagerup
property (see \cite{CJV}).

\end{section}

\begin{section}{The amalgamated case}
 In this short section, we briefly explain how to extend the previous
results to amalgamated free products.

 We start with some considerations about the
amalgamated Haagerup tensor product
and Hilbert $C^{*}$-modules. A good reference for this material is \cite{BLM}.

For simplicity, we will always assume that all $C^*$-algebras are
unital  and all modules are full. If $H$ and $K$ are two right
$C^*$-modules over a $C^*$-algebra $B$, $\B_B(H, K)$ denotes the
space of all adjointable $B$-right-modular maps from $H$ to $K$,
and $\K_B(H, K)$ is its subspace of compact maps.

 Let $H$ be a right $C^*$-$B$-module. Then there is a natural antilinear
map $t:H\to \K_B(H,B)$ defined by $t_{h}(x)=\la h,\;x\ra$. Using
this identification, $h \leftrightarrow t_h$, $H$ inherits a left
$C^*$-$B$-module from that of $\K_B(H,B)$. For instance, one has
$b\cdot t_h=t_{hb^*}$, and $\la t_h,\; t_{h'}\ra= \la
h,\;h'\ra_{H}$. We will omit the $t$ in this identification and we
call this module $\overline H$. Similarly, if $K$ is a left
$C^*$-$B$-module, we can define an associated right
$C^*$-$B$-module $\overline K$.

Let $H$ be a right $C^*$-$B$-module. The $C^*$-module structure
can be used to endow $H$ with a canonical operator space
structure, just as the column Hilbert space structure for a
Hilbert space. More precisely, it is inherited from the
identifications
 $$H=\K_B(B,H), \quad \quad  \forall n\ge 1,\; \M_n(H)
 =\K_B(B^n,H^n).$$
Here $B^n$ (resp. $H^n$) is the direct sum of $n$ copies of $B$
(resp. $H$). For instance, the scalar product of $B^n$  is
$\langle (b_i), (b'_i)\rangle=\sum_{i=1}^n b_i^*b^{\prime}_i$.
(Note that $\K_B(B,H)=\B_B(B,H)$.)  We will denote this operator
space structure by $H_C$.

Similarly, if $K$ is a left $C^*$-$B$-module,  the row operator
space structure on $K$ comes from the identifications
 $$K=\K_B(\overline{K},B), \quad \quad  \forall n\ge 1,\;
 \M_n(K)=\K_B(\overline{K}^n, B^n).$$
We will denote this operator space structure by
$K_R$.

 The Haagerup tensor product has a modular counterpart for operator spaces
$X$ and $Y$ which are respectively right and left $B$-modules ($B$ can merely
be an algebra). It is defined as the quotient of the usual Haagerup product
by the subspace spanned by $xb\tens y-x\tens by$ with $x\in X$, $y\in Y$ and
$b\in B$, we denote it by $X\tens_{hB}Y$.
 Obviously, it has the universal property that for any completely contractive
maps $\sigma:X\to \B(H)$, $\rho:Y\to \B(H)$ (with $H$ a Hilbert
space) so that $\sigma(xb)\rho(y)=\sigma(x)\rho(by)$ (with the
above notations), there is a unique completely contractive map
$\sigma.\rho:X\tens_{hB}Y\to \B(H)$ so that $\sigma.\rho(x\tens
y)=\sigma(x)\rho(y)$.

\medskip

 Now, assume that $B$ is a unital $C^*$-algebra, and
$(A_i)_{i\in I}$ is a family of $C^*$-algebras containing $B$ as a
unital $C^*$-subalgebra with a conditional expectation
$\phi_i:A_i\to B$ whose GNS representation is faithful. Then we
can define the reduced free product of the $A_i$ with amalgamation
over $B$. The space $H_i=L_2(A_i,\phi_i)$, (resp. $H_i^{op}$) is
the right Hilbert (resp. left) $C^*$-$B$-module  obtained from $A$
using the product
 $$\la\hat a,\;\hat b\ra_{ H_i}=\phi_i(a^*b)\qquad (\textrm{resp.}
 \la\hat a,\;\hat b\ra_{ H_i^{op}}=\phi_i(ab^*)).$$
As before there
is a natural identification of $H_i^{op}$ with $\K_B(H_i,B)$ using
the duality
 $$\la\hat a,\;\hat b\ra_{(H_i^{op},\; H_i)}=\phi_i(ab).$$
Equivalently, this means that $H_i^{op}$ is isometric to
$\overline{H_i}$.
 Actually, $H_i$ (resp. $H_i^{op}$) is also a left (resp. right)
$A_i$-module with the natural multiplication. To be consistent
with our previous notations in the non amalgamated case, we denote
by $\m A_i$ and $\m H_i$ the kernel of $\phi_i$ on $A_i$ and
$H_i$.

Let $\F$ be the Hilbert $C^*$-$B$-module Fock space associated to
the free product:
 $$\F=\mathbb{C}\cdot\Omega ~~\bigoplus~\sd n \plB H i n,$$
where $\tens_B$ stands for the interior tensor product with
respect to the left action of $B$.

 The reduced amalgamated free product
$*_{i\in I}(A_i,\phi_i)$ is then the $C^*$-subalgebra of
$\B_B(\F)$ generated by the copies of the $A_i$'s (cf. section 1
for the formulas). In the algebraic amalgamated free product of
the $A_i$, we will use $\Sigma_d$ to denote the space of
homogeneous elements of length $d$.

In this setting, the projections $P_k$ of section 2 are well
defined as adjointable right $B$-modular maps on $\F$. Moreover,
the $P_k$ commute with the action of $B$.  We keep all the
previous notations, $L_1$ is the operator space in $\B_B(\F)$
spanned by $(P_k \m A_k P^\bot_k)_{k\in I}$, and $K_1$ that
spanned by $(P_k^\bot \m A_k P_k)_{k\in I}$. $L_1$ and $K_1$ are
naturally $B$-bimodules. Lemmas \ref{Colonne} and \ref{Ligne}
remain true with the same proof. Also, the space
$\ell_\infty((A_i))$ is naturally a $B$-bimodule and lemma
\ref{linfi} still holds.

Then we introduce an operator space structure on $
\Sigma_1^{\tens^d}$ via the following inclusion
 $$\iota : \left\{\begin{array}{ccccc}
 \Sigma_1^{\tens^d} & \to &
 \bigoplus_{k=0}^d L_1^k\tens_{hB} K_1^{d-k} & \bigoplus_\infty &
 \bigoplus_{k=0}^{d-1}L_1^k\tens_{hB}\ell_\infty((A_i))
 \tens_{hB} K_1^{d-k-1}\\
 a& \mapsto& \Big((\theta_1^k\tens\rho_1^{d-k}(a))_{k=0}^d&, &
 (\theta_1^k\tens \Id\tens \rho_1^{d-k-1}(a))_{k=0}^{d-1}\Big)
 \end{array} \right.,$$
where $L_1^k=L_1^{\tens_{hB}^ k}$ and $K_1^k=K_1^{\tens_{hB}^ k}$.

 From the extension of the lemmas in section 1 and the universal property
of the modular Haagerup tensor product, it is obvious that the
majoration of Theorem \ref{khintchine} remains true. To adapt the
proof of the minoration we simply need to show the following,
which is just the modular version of the well-known fact that
$H_C\tens_hA\tens_h K_R= A\tens_{\min}\K(\overline{K},H)$.

\begin{prop}
Let $H,\, K$ and $X$ be right Hilbert $C^*$-$B$-modules. Assume
that $X$ is also a left $B$-module and $A\subset \B_B(X)$ is a
closed $B$-bimodule (i.e. stable under left and right
multiplications by elements from $B$). Then the map
  $$ \Phi: \left\{ \begin{array}{ccc}
  H_C\tens_{hB} A\tens_{hB} \overline{K}_R &
  \to & \B_B(K\tens_B X, H\tens_B X) \\
  h\tens a\tens\bar k   &\mapsto &
  \big( k'\tens x \mapsto h\tens (a \,\la k,\;k'\ra)x\big)
  \end{array}\right.$$
is a complete isometry.
\end{prop}

\begin{pf}
 From the universal property
of Haagerup tensor product, $\Phi$ is completely contractive. Let
$x=\sum_i h_i\tens a_i\tens \bar k_i\in H_C\tens A\tens
\overline{K}_R$. It is know that there are contractive
approximate units $(e_n)$ in $\K_B(H)$ and $(f_m)$ in $\K_B(K)$,
which are of the following form :
 $$e_n=\sum_{s=1}^{M_n} u_s\tens u_s, \qquad
  f_n=\sum_{t=1}^{N_n} v_t\tens v_t$$
(see \cite[8.1.23]{BLM}). Let $T_{m,n}(x)=\sum_{i} e_m(h_i)\tens
a_i\tens \bar f_n(\bar k_i)$. It is clear that $\lim_{m,n}
T_{m,n}(x)=x$ in the Haagerup tensor product, so one only needs to
check that $\|T_{m,n}(x)\|_h\le \|\Phi(x)\|$. To this end we write
 \begin{eqnarray*}
 T_{m,n}(x)&=&\sum_{i} e_m(h_i)\tens a_i\tens \bar f_n(\bar k_i) \\
 &=& \sum_{i,s,t} u_s \tens \la g_{s},\;h_i\ra\, a_i\,
 \la k_i,\;v_{t}\ra  \tens  \bar v_t\\
 &=& \sum_{s,t} u_s \tens \sum_i\la  u_s,\;h_i\ra\,a_i\,
  \la k_i,\;v_{t}\ra \tens \bar v_t
 \end{eqnarray*}
  Thus $T_{m,n}(x)$ can be written as $UWV$, with $U$ a row matrix corresponding
to $(u_s)_{s}$, $V$ a column matrix with entries $(\bar v_t)_{t}$
 and $W$ the  matrix
$\big(\sum_i \la g_s,\;h_i\ra \,a_i\, \la k_i,\;
v_t\ra\big)_{s,t}$. The matrices $U$ and $V$ have norm less than
one as the approximate units are contractive.
 To get the conclusion we need to prove that $\|W\|\le
\|\Phi(x)\|$, where $W$ is viewed as an operator from $X^{N_n}$ to
$X^{M_m}$. For elements  $\xi=(\xi_t)\in  X^{N_n}$ and
$(\eta_s)\in X^{M_m}$, we have
 \begin{eqnarray*}
 \la \eta,\; W\xi\ra&=& \sum_{s,t,i} \la  \eta_s,\;
 \la  u_s,\;h_i\ra\, a_i\, \la  k_i,\;v_t\ra\, \xi_t \ra_{X}\\
 & = & \la \sum_s u_s \tens \eta_s,\; \Phi(x) \sum_t v_t\tens
 \xi_t\ra_{H\tens_BX}\,.
 \end{eqnarray*}
The conclusion then follows easily from the fact that both $(\la
u_s,\;u_{s'}\ra)_{s,s'\le M_m} \in \M_{M_m}(B)$ and  $(\la
v_t,\;v_{t'}\ra)_{t,t'\le N_n} \in \M_{N_n}(B)$ are contractive.
\end{pf}

>From this result, it follows that Khintchine inequalities are true in the
amalgamated case.

\medskip

Now let $A$ and $B$ be $C^*$-algebras with $B\subset A$ a unital
subalgebra. Assume that there is a conditional expectation $\phi$
from $A$ onto $B$. Then $\phi$ induces a right Hilbert
$C^*$-$B$-module in $A$. $A$ is, of course,  a left $B$-module in
the natural way. We also view $A$ as a subalgebra of $\B_B(A)$.

\begin{prop}
Let $A$, $B$ be as above, and let $H, K$ be right
$C^*$-$B$-modules. If $A$ is exact, then  $H_C\tens_{hB}
A\tens_{hB} \overline{K}_R$ is 1-exact.
\end{prop}

\begin{pf}
 Since exactness is a local notion and the Haagerup tensor product
is injective, we can assume that $H$ and $K$ are countably
generated. Hence by Kasparov's stabilization theorem, they are
complemented in $B^\N$. Then
 $$H_C\tens_{hB} A\tens_{hB} \overline{K}_R\subset
 B^\N_C\tens_{hB} A\tens_{hB}\overline{B^\N}_R\quad\textrm{
 completely isometric}.$$
 Thus it suffices to show the result for
$H=K=B^\N$. In the same way, exactness is stable by inductive
limit, by the same kind of arguments we are reduced to show it for
$H=K=B^n$ ($n\ge 1$). But then the previous proposition gives that
the Haagerup tensor product is exactly $\M_n(A)$, which is 1-exact.
\end{pf}

\medskip

Using this proposition and the Khintchine inequalities, we
conclude, as in section~4, that  the exactness is stable under
amalgamated reduced free products.

\end{section}

\end{document}